\documentclass[11pt]{amsart}

\usepackage{amsmath}
\usepackage{amsthm}
\usepackage{amssymb}
\usepackage{upref}
\usepackage{amscd}
\usepackage[all]{xy}

\def\Q{\mathbb Q\mspace{1mu}}
\def\Z{\mathbb Z}
\def\F{\mathbb F}

\def\C{\mathbb C}
\def\R{\mathbb R}

\def\H{\mathbb H}

\def\Qa{\overline{\mathbb Q}\mspace{1mu}}
\def\Qp{{\mathbb Q}_p}
\def\Qpa{\overline{\mathbb Q}_p}
\def\Fpa{\overline{\mathbb F}_p}

\def\l{\ell}
\def\DcrK{{\rm\bf D}^*_{\text{cris},K}}

\DeclareMathOperator{\ord}{ord}
\DeclareMathOperator{\lcm}{lcm}
\DeclareMathOperator{\id}{Id}
\DeclareMathOperator{\im}{Im}
\DeclareMathOperator{\Ker}{Ker}
\DeclareMathOperator{\tr}{Tr}
\DeclareMathOperator{\Endom}{End}
\DeclareMathOperator{\Hom}{Hom}
\DeclareMathOperator{\aut}{Aut}
\DeclareMathOperator{\ist}{Ist}

\DeclareMathOperator{\pchar}{P_{char}}
\DeclareMathOperator{\pmin}{P_{min}}

\DeclareMathOperator{\prd}{Prd}
\DeclareMathOperator{\roots}{Roots}
\DeclareMathOperator{\diag}{Diag}
\DeclareMathOperator{\spec}{Spec}
\DeclareMathOperator{\ind}{ind}
\DeclareMathOperator{\gal}{Gal}
\DeclareMathOperator{\Frac}{Frac}
\DeclareMathOperator{\inv}{inv}
\DeclareMathOperator{\pol}{Pol}

\DeclareMathOperator{\M}{{\bf M}}
\DeclareMathOperator{\W}{{\bf W}}
\DeclareMathOperator{\D}{{\bf D}}

\DeclareMathOperator{\Fil}{Fil}

\begingroup
\theoremstyle{plain}
\newtheorem{thm}{Theorem}[section]
\newtheorem{prop}[thm]{Proposition}
\newtheorem{lem}[thm]{Lemma}
\newtheorem{cor}[thm]{Corollary}

\theoremstyle{definition}
\newtheorem{defn}[thm]{Definition}

\newtheorem{ex}[thm]{Example}
\newtheorem{rem}[thm]{Remark}

\endgroup

\newenvironment{mthm}{\medskip\noindent{\bf Theorem.}}{\medskip}
\newenvironment{ack}{\medskip\noindent{\bf Acknowledgements.}}{\medskip}

\begin{document}

\title[A class of representations arising from abelian varieties over $\Qp$]{A class of $p$-adic Galois representations \\
arising from abelian varieties over $\Qp$}

\author{Maja Volkov}

\begin{abstract}
Let $V$ be a $p$-adic representation of the absolute Galois group $G$ of $\Qp$ that becomes crystalline over a finite tame extension, and assume $p\neq 2$. We provide necessary and sufficient conditions for $V$ to be isomorphic to the $p$-adic Tate module $V_p({\mathcal A})$ of an abelian variety ${\mathcal A}$ defined over $\Qp$. These conditions are stated on the filtered $(\varphi,G)$-module attached to $V$. 
\end{abstract}

\curraddr{18 rue Dussoubs, 75002 Paris, France}

\email{volkov@math.uni-muenster.de}

\maketitle

\medskip

\tableofcontents

\section*{Introduction}
Fix a prime number $p$ and an algebraic closure $\Qpa$ of $\Qp$. Let ${\mathcal A}$ be an abelian variety over $\Qp$ of dimension $d$. Let ${\mathcal A}[p^n]$ be the group of $p^n$-torsion points with values in $\Qpa$ and $T_p({\mathcal A})=\varprojlim{\mathcal A}[p^n]$ the $p$-adic Tate module of ${\mathcal A}$: it is a free $\Z_p$-module of rank $2d$ on which the Galois group $G=\gal(\Qpa /\Qp)$ acts linearly and continuously. The $p$-adic representation of $G$ attached to ${\mathcal A}$ (also called Tate module) is 
$$V_p({\mathcal A}) = \Qp\otimes_{\Z_p}T_p({\mathcal A}).$$
A $p$-adic representation $V$ of $G$ arises from an abelian variety over $\Qp$ if there exists an ${\mathcal A}/\Qp$ such that $V\simeq V_p({\mathcal A})$ as $G$-modules. Such representations are classical objects known to have motivic properties. We want to consider an inverse problem, namely, to detect which $p$-adic representations $V$ of $G$ arise from abelian varieties. In this paper we solve this problem for the class of representations that become crystalline over a tame extension of $\Qp$, and under the assumption $p\neq 2$. The solution is provided in terms of necessary and sufficient conditions on the filtered $(\varphi,G)$-module attached to $V$.   

\medskip
Being crystalline is the $p$-adic analogue of the notion of good reduction for $\l$-adic representations, $\l\neq p$. Fontaine constructed a functor establishing an equivalence of categories between potentially crystalline $p$-adic representations of $G$ and weakly admissible filtered $(\varphi,G)$-modules. The latter are finite-dimensional vector spaces equipped with a semilinear Frobenius map $\varphi$, a semilinear action of a finite quotient of $G$ commuting with $\varphi$, and a $G$-stable filtration. Weak admissibility is a relation between $\varphi$ and the filtration ensuring that the object comes from an actual representation of $G$. Forgetting the filtration we obtain a $(\varphi,G)$-module to which we can associate functorially a representation of the Weil group $W={\rm Weil}(\Qpa /\Qp)$, in such a way that isomorphism classes are preserved.

\medskip
The abelian variety ${\mathcal A}/\Qp$ has potential good reduction if and only if $V_p({\mathcal A})$ is potentially crystalline. Hence the above theory applies, and the varieties involved here have potential good reduction. It should be pointed out that we are using contravariant functors. Since $V_p({\mathcal A})$ is dual to $H^1_{\text{\'et}}({\mathcal A}_{\Qpa},\Qp)$ these functors are covariant on the $H^1_{\text{\'et}}$'s. Let ${\mathcal A}$ have good reduction over the finite tame Galois extension $K/\Qp$ with residue field $k$ and maximal unramified subfield $K_0$. Let $(D,\Fil)$ be the filtered $\bigl(\varphi,\gal(K/\Qp)\bigr)$-module attached to $V=V_p({\mathcal A})$ and $\Delta$ the Weil representation attached to $D$. Then $(D,\Fil)$ satisfies the following conditions:  
\begin{itemize}
\item[(1)] A lifting of the geometric Frobenius acts semisimply on $\Delta$ and its characteristic polynomial is a $p$-Weil polynomial (see definition~\ref{pWpoldefn})
\item[(2)] $\Delta$ is defined over $\Q$ and has Tate type (see definition~\ref{Tatetypedefn})
\item[(3)] There exists a nondegenerate skew form $(D,\Fil)\times (D,\Fil)\rightarrow K_0\{-1\}$
\item[(4)] $(D,\Fil)$ has Hodge-Tate type $(0,1)$.
\end{itemize}
Write $A$ for the special fibre of the N\'eron model of ${\mathcal A}_{K}$. Conditions $(1)$, $(2)$, and $(3)$ unfiltered express the fact that the Weil representation $W\rightarrow \aut_{K_0}(\Delta)$ factors through a map $W\rightarrow \bigl(\Q\otimes_{\Z}\Endom_{k}(A)\bigr)^{\times}$. In particular the descent datum furnished by the action of $\gal(K/\Qp)$ on $(D,\Fil)$ is ``geometric''. More precisely $(1)$ follows from the Weil conjectures for abelian varieties over $\F_p$ and $(2)$ from the results of Serre and Tate on potential good reduction of abelian varieties together with Tate's description of the algebra $\Q\otimes_{\Z}\Endom_{k}(A)$. Having Tate type means that the dimension of some distinguished subobjects of $\Delta$ is divisible by an arithmetic invariant of these. Condition $(3)$ means that $V$ may be equipped with a $G$-equivariant symplectic form $V\times V\rightarrow \Qp(1)$, and indeed any polarisation on ${\mathcal A}$ induces such a form. Finally $(4)$ is a well-known condition on the filtration. Now our result is that these necessary conditions turn out to be sufficient.

\begin{mthm} (see theorem~\ref{mainthm}) {\em Let $p\neq 2$. Let $V$ be a $p$-adic representation of $G$ that becomes crystalline over a finite tame extension. The following are equivalent:
\begin{itemize}
\item[(i)] There exists an abelian variety ${\mathcal A}$ over $\Qp$ such that $V\simeq V_p({\mathcal A})$
\item[(ii)] The filtered $(\varphi,G)$-module attached to $V$ satisfies $(1)$, $(2)$, $(3)$, and $(4)$.
\end{itemize}}
\end{mthm}

Given a $p$-adic representation $V$ of $G$ satisfying the conditions of the theorem we need to construct an abelian variety ${\mathcal A}$ over $\Qp$ such that $V\simeq V_p({\mathcal A})$. Let $V$ be crystalline over the tame Galois extension $K/\Qp$. We first construct the relevant objects defined over finite fields. The Honda-Tate theory enables us to produce from $(1)$ an abelian variety $A_0$ over $\F_p$ having the right Frobenius. Using Tate's theorems we produce essentially from $(2)$ an appropriate automorphism of $A=A_0\times_{\F_p}k$ furnishing a geometric descent datum. We then construct the relevant objects in characteristic zero. By a result of Breuil (which is why we assume $p$ odd) condition $(4)$ furnishes a $p$-divisible group over the ring of integers $O_K$ of $K$ lifting $A(p)$. According to the Serre-Tate theory of liftings this produces a formal abelian scheme over $O_K$. Then condition $(3)$ enables us to lift some polarisation on $A$, which insures by Grothendieck's criterion that we obtain an algebraic scheme, hence an abelian variety over $K$ with good reduction. Finally a Galois descent criterion, applicable thanks to the geometricity of the descent datum, shows that this abelian variety is defined over $\Qp$ and has the right Tate module.

\medskip
The paper is organised as follows. Section~\ref{sec:pgr} contains assorted facts on potential good reduction. We review the needed representation theory in~\ref{sec:repn} and potential good reduction of abelian varieties in~\ref{sec:inertia}. We recall the Galois descent criterion in~\ref{sec:desc} and the influence of potential good reduction on duals and polarisations in~\ref{sec:dual}. Section~\ref{sec:abvarfinite} deals with the Honda-Tate theory. It is reviewed in~\ref{sec:HT} and in~\ref{sec:HTalg} we prove a technical result on endomorphism algebras (proposition~\ref{embedprop}). In section~\ref{sec:galpairs} we introduce the appropriate objects yielding representations with geometric descent: they consist of an abelian variety over $\F_p$ together with a finite subgroup of $\Fpa$-automorphisms, subject to a condition of Galois nature. These objects, that we call Galois pairs, are defined in~\ref{sec:defn} (definition~\ref{galabdefn}), their associated representations are constructed in~\ref{sec:repnsgalpairs}, and a decomposition result (proposition~\ref{decompprop}) is proved in~\ref{sec:tamegalpairs}. In section~\ref{sec:repnsabvarFp} we describe such representations in the tame case. Some arithmetic properties of semisimple representations defined over $\Q$ are proved in~\ref{sec:ssQ}. The determination of tame representations arising from Galois pairs is carried out in~\ref{sec:galpairsrepns} (theorem~\ref{HTthm}), followed by some examples in~\ref{sec:ex}. Section~\ref{sec:repabvar} deals with the lifting procedure. We review in~\ref{sec:pol} the needed facts on polarisations and Rosati involutions, and prove in~\ref{sec:polgalpairs} the existence of suitable polarisations on tame Galois pairs (proposition~\ref{polgalprop}). In~\ref{sec:geomsp} we show that such polarisations may be lifted when the representation associated to a Galois pair lifts to a symplectic $G$-module (proposition~\ref{polprop}). We then wrap everything up to prove the main theorem in~\ref{sec:thm} (theorem~\ref{mainthm}).

\begin{ack}
It is a pleasure to thank A.J.~Scholl, J.~Nekov\'a\v r, N.~Naumann, and M.~Kisin for helpful conversations. Special thanks to O.~B\"ultel for lemma~\ref{oliverlem} and to E.~Frossard for lemma~\ref{sbalglem}. 
\end{ack}

\section{Potential good reduction}
\label{sec:pgr}

\subsection{Representations}
\label{sec:repn}
Let $I$ be the inertia subgroup of $G=\gal(\Qpa /\Qp)$ and $W={\rm Weil}(\Qpa /\Qp)$ the Weil group. For extensions $K$ of $\Qp$ (assumed to be contained in $\Qpa$) write $G_K=\gal(\Qpa/K)$, $I_K=I(\Qpa/K)$, and $W_K={\rm Weil}(\Qpa/K)$. Let $\Qp^{\text{un}}$ be the maximal unramified extension of $\Qp$, $\Fpa$ its residue field, $\sigma$ the absolute Frobenius, and $W(k)$ the ring of Witt vectors with coefficients in $k\subseteq\Fpa$. We also fix an algebraic closure $\Qa$ of $\Q$ and a field embedding $\Qa \hookrightarrow \Qpa$.

\subsubsection{Filtered $(\varphi,G)$-modules}
\label{sec:filmod}
Let $K$ be a finite Galois extension of $\Qp$ and $K_0$ the maximal unramified subfield $K$. A filtered $\bigl(\varphi,\gal(K/\Qp)\bigr)$-module $(D,\Fil)$ is a finite dimensional $K_0$-vector space $D$ equipped with
\begin{itemize}
\item[{\scriptsize $\bullet$}] a $\sigma$-semilinear bijective Frobenius map $\varphi : D \xrightarrow{\sim} D$ 
\item[{\scriptsize $\bullet$}] a $\sigma$-semilinear action of $\gal(K/\Qp)$ commuting with $\varphi$
\item[{\scriptsize $\bullet$}] a decreasing filtration $\Fil=(\Fil^iD_K)_{i\in \Z}$ on $D_K=K\otimes_{K_0}D$ by $\gal(K/\Qp)$-stable subspaces such that $\Fil^iD_K=D_K$ for $i\ll 0$ and $\Fil^iD_K=0$ for $i\gg 0$.
\end{itemize}
Such objects form a category (\cite{Fo2}), the morphisms being $K_0$-linear maps commuting with the action of $\varphi$ and $\gal(K/\Qp)$, and preserving the filtration after scalar extension to $K$. The dual of $(D,\Fil)$ is the $K_0$-linear dual $D^*$ with $\varphi^*f=\sigma f\varphi^{-1}$ and $g^*f=gfg^{-1}$ for $f\in D^*$, $g\in\gal(K/\Qp)$, and $\Fil^iD_K^*$ consists of linear forms on $D_K$ vanishing on $\Fil^jD_K$ for all $j>-i$. The Tate twist $D\{-1\}$ of $(D,\Fil)$ is $D$ as a $K_0$-vector space with the same action of $\gal(K/\Qp)$ and $\varphi\{-1\}=p\varphi$, $\Fil^i(D\{-1\})_K=\Fil^{i-1}D_K$. The object $(D,\Fil)$ has Hodge-Tate type $(0,1)$ if the filtration jumps in degrees $0$ and $1$, meaning $\Fil^iD_K=D_K$ if $i\leq 0$, $\Fil^iD_K=0$ if $i\geq 2$, and $\Fil^1$ is a nontrivial subspace.

\medskip
Now let $V$ be a $p$-adic representation of $G$. Put 
$$\DcrK(V) \underset{\text{def}}{=} \Hom_{\Qp[G_K]}(V,B_{\text{cris}}).$$
Then $\DcrK(V)=(D,\Fil)$ is a filtered $\bigl(\varphi,\gal(K/\Qp)\bigr)$-module (\cite{Fo1}, \cite{Fo2}). We always have $\dim_{K_0}D\leq \dim_{\Qp}V$ and $V$ is crystalline over $K$ when equality holds. In this case $(D,\Fil)$ satisfies a property relating $\varphi$ with the filtration, called weak admissibility. The association $V\mapsto \DcrK(V)$ is functorial (contravariant) and sets an anti-equivalence of tannakian categories between $p$-adic representations of $G$ crystalline over $K$ and weakly admissible filtered $\bigl(\varphi,\gal(K/\Qp)\bigr)$-modules (\cite{Fo2}, \cite{Co-Fo}). We have $\DcrK\bigl(V(1)\bigr)=\DcrK(V)\{-1\}$ where $V(1)=\Qp(1)\otimes_{\Qp}V$ is the usual Tate twist.

A representation is potentially crystalline if it is crystalline over some finite extension. Being crystalline depends only on the action of inertia (\cite{Fo2} 5.1.5). Let $K\subseteq F$ be two finite Galois extensions such that $F/K$ is unramified, with respective maximal unramified subfields $K_0\subseteq F_0$. Then $V$ is crystalline over $K$ if and only if it is crystalline over $F$. Let $(D,\Fil)$ be a filtered $\bigl(\varphi,\gal(F/\Qp)\bigr)$-module. Then $(D,\Fil)=(F_0\otimes_{K_0}D_0,F\otimes_{K}\Fil_0)$ where $(D_0,\Fil_0)$ is the filtered $\bigl(\varphi,\gal(K/\Qp)\bigr)$-module obtained from $(D,\Fil)$ by taking $\gal(F/K)$-invariants.

\medskip
Let ${\mathcal A}_0$ be an abelian variety over $\Qp$ having potential good reduction. Then $V_p({\mathcal A}_0)$ is potentially crystalline, and it is crystalline over $K$ if and only if ${\mathcal A}_0$ has good reduction over $K$ (\cite{Co-Io}, see also~\cite{Br}). Then $\DcrK\bigl(V_p({\mathcal A}_0)\bigr)=(D,\Fil)$ has Hodge-Tate type $(0,1)$ and $\dim_K\Fil^1D_K=\dim{\mathcal A}_0$.

\subsubsection{Dieudonn\'e modules}
\label{sec:dieumod}
Let $(D,\Fil)$ be a filtered $\bigl(\varphi,\gal(K/\Qp)\bigr)$-module. Forgetting the filtration leads to the obvious notion of a $\bigl(\varphi,\gal(K/\Qp)\bigr)$-module $D$, and forgetting in addition the action of $\gal(K/\Qp)$ to a $\varphi$-module. Let $k$ be the residue field of $K$. The contravariant functor $\Gamma \mapsto \D(\Gamma) = K_0\otimes_{W(k)}\M(\Gamma)$ described in~\cite{Fo4} sets an anti-equivalence between the category of $p$-divisible groups over $k$ up to isogeny and $\varphi$-modules with slopes between $0$ and $1$.

\medskip
Let ${\mathcal A}_0\times_{\Qp}K$ have good reduction, let $A$ be the special fibre of the N\'eron model and $A(p)$ its associated $p$-divisible group over $k$. Put
$$D(A) \underset{\text{def}}{=} \D\bigl(A(p)\bigr).$$
If $\DcrK\bigl(V_p({\mathcal A}_0)\bigr)=(D,\Fil)$ there is a canonical isomorphism of $\varphi$-modules $D(A)\simeq D$ allowing us to identify these two objects (\cite{Fo5}). The Frobenius $\varphi$ on $D(A)$ is induced by the Frobenius endomorphism of $A$. The filtration $\Fil$ is a lifting datum, that is, it carries the information that $A(p)$ is the special fibre of a $p$-divisible group over the ring of integers of $K$, namely the one of the N\'eron model. The action of $\gal(K/\Qp)$ is a descent datum reflecting the fact that we actually started with an object defined over $\Qp$.

\subsubsection{Weil representations}
\label{sec:Wrepns}
A $\bigl(\varphi,\gal(K/\Qp)\bigr)$-module is a semilinear object (unless $K_0=\Qp$) and it is convenient have a linear one at hand naturally attached to it. The Weil group sits in a short exact sequence
$$1 \longrightarrow I \longrightarrow W \overset{\upsilon}{\longrightarrow} \Z \longrightarrow 1$$
where $\upsilon$ sends a lifting of the arithmetic Frobenius to $1$. From a given $\bigl(\varphi,\gal(K/\Qp)\bigr)$-module $D$ we get a $K_0$-linear representation $\W(D)=(\Delta,\rho)$ of $W$ as follows (\cite{Fo3}). As a $K_0$-vector space $\Delta$ is $D$, and the action is given by 
$$\rho(w)= (w \bmod W_K) \varphi^{-\upsilon(w)} \in \aut_{K_0}(\Delta), \quad w\in W.$$
We have $I_K\subseteq \Ker\rho$, that is, $(\Delta,\rho)$ has good reduction over $K$. It is tame if $I_p\subset \Ker\rho$ with $I_p$ the maximal pro-$p$-group contained in $I$. The association $D\mapsto \W(D)$ is functorial (covariant), and we have a canonical isomorphism 
$$K_0\otimes_{\Qp}\Hom_{\varphi}(D_1,D_2) \simeq \Hom_{_{W}}\bigl(\W(D_1),\W(D_2)\bigr)$$
so that isomorphism classes are in bijection. We may always assume that $K$ is a Galois extension of minimal degree over which $\Delta$ has good reduction so that $K_0$ is minimal as well (\cite{Fo3} 1.3.7).

\medskip
A semisimple $(\Delta,\rho)$ is defined over $\Q$ if $\tr\bigl(\rho(w)\bigr)\in \Q$ for all $w\in W$. In this case it can be compared with $\Q_{\l}$-linear Weil representations with $\l\neq p$ (\cite{Fo3}). For each prime $\l\neq p$ let $\Delta_{\l}$ be a $\Q_{\l}$-linear semisimple Weil representation defined over $\Q$ (same definition) and let $\Delta_p$ be a $K_0$-linear one. Then the system $(\Delta_{\l})_{\l}$ is compatible if the trace maps are independent of $\l$. It is known that abelian varieties with potential good reduction give rise to such systems.

\subsection{The action of inertia}
\label{sec:inertia}
Let ${\mathcal A}_0$ be an abelian variety over $\Qp$ having potential good reduction. Recall (\cite{Se-Ta}, \cite{Gr2}) that there is a smallest extension $M$ of $\Qp^{\text{un}}$ over which ${\mathcal A}_0$ acquires good reduction. It is Galois, given by 
$$G_M = I_M = I \cap \Ker(\rho_{\l})$$ 
where $\l$ is any prime different from $p$ and $\rho_{\l} : G \rightarrow \aut_{\Q_{\l}}\bigl( V_{\l}({\mathcal A}_0)\bigr)$ the $\l$-adic representation. Since $M$ is Galois over $\Qp$ it can be descended to a finite totally ramified extension $L$ of $\Qp$, that is, $L\Qp^{\text{un}} = M$. The Galois closure $K$ of $L$ enjoys the following properties: $K/L$ is unramified, $K\Qp^{\text{un}}=M$, $I_K=I_M$, and $\gal(K/\Qp) \simeq I(K/\Qp)\rtimes \gal(k/\F_p)$ (semidirect product) with $k$ the residue field of $K$. Furthermore ${\mathcal A}_0$ has good reduction over $L$ and $K$ is a Galois extension of minimal degree over which ${\mathcal A}_0$ acquires good reduction.

\medskip
The group $\gal(K/\Qp)$ acts on ${\mathcal A}={\mathcal A}_0\times_{\Qp}K$ via its action on $K$. This extends functorially to an action on the N\'eron model of ${\mathcal A}$, and on its special fibre $A/k$. The inertia subgroup acting trivially on $k$ it acts on $A$ by $k$-automorphisms, i.e. by an antimorphism 
$$\nu : I(K/\Qp) \rightarrow \aut_k(A).$$ 
This map is injective because of the minimality of $K$. Indeed, $V_p({\mathcal A}_0)$ is crystalline over $K^{\Ker(\nu)}$ so ${\mathcal A}_0$ has good reduction over this field. Hence $\nu$ identifies $I(K/\Qp)$ with a finite subgroup of $\aut_k(A)$.

\medskip
Of course this is nothing else than the classical Weil criterion (\cite{We}). For each $g\in \gal(K/\Qp)$ let ${\mathcal A}^g$ be the twisted abelian variety over $K$, with the usual relation ${\mathcal A}^{gh}=({\mathcal A}^h)^g$. Since ${\mathcal A}$ is defined over $\Qp$ there exists a system of $K$-isomorphisms
$$f_g : {\mathcal A} \xrightarrow{\sim} {\mathcal A}^g , \quad g \in \gal(K/\Qp)$$
satisfying the cocycle condition
$$f_{gh} = (f_h)^g \circ f_g \quad \text{ for all } g,h \in \gal(K/\Qp).$$
Now if $g \in I(K/\Qp)$ the special fibre of ${\mathcal A}^g$ is $A$, so by restriction to the inertia subgroup and passing to special fibres we recover $\nu : I(K/\Qp) \hookrightarrow \aut_k(A)$. 

\begin{rem}
Since ${\mathcal A}_0$ acquires good reduction over the totally ramified extension $L$ we have $A=A_0\times_{\F_p}k$, where $A_0/\F_p$ is the special fibre of the N\'eron model of ${\mathcal A}_0 \times_{\Qp}L$.
\end{rem}

Thus $\DcrK\bigl(V_p({\mathcal A}_0)\bigr)=(D,\Fil)$ is a $\bigl(\varphi,\gal(K/\Qp)\bigr)$-filtered module with an action of $I(K/\Qp)$ preserving the filtration and coming from the composite morphism
$$I(K/\Qp) \underset{\nu}{\hookrightarrow} \aut_k(A) \underset{\text{can}}\hookrightarrow \aut_{\varphi}(D).$$

\subsection{Galois descent}
\label{sec:desc}
The following result is a converse to the situation explained in section~\ref{sec:inertia}. Let $F$ be a finite extension of $\Qp$. 

\begin{thm}[\cite{Vo} Thm.4.5]
\label{descthm}
Let $K/F$ be a totally ramified finite Galois extension with residue field $k$. Let ${\mathcal A}/K$ be an abelian variety having good reduction with special fibre $A/k$. Then ${\mathcal A}$ is defined over $F$ if and only if the action of $G_K$ on $T_p({\mathcal A})$ extends to an action of $G_F$ coming from $k$-automorphisms of $A$, that is, from a morphism
$$\nu : \gal(K/F)=I(K/F) \rightarrow  \aut_k(A).$$
\end{thm}

Specifically, there is an abelian variety ${\mathcal A}_0/F$ and a $K$-isomorphism $\psi: {\mathcal A}_0 \times_{F}K \xrightarrow{\sim} {\mathcal A}$ inducing a $G_F$-equivariant isomorphism $T_p(\psi): T_p({\mathcal A}_0) \xrightarrow{\sim} T_p({\mathcal A})$, with $G_F$ acting naturally on $T_p({\mathcal A}_0)$ and by the given extended action on $T_p({\mathcal A})$. The pair $({\mathcal A}_0,\psi)$ is unique up to $F$-isomorphism. Also, the smallest extension of $F$ contained in $K$ over which ${\mathcal A}_0$ acquires good reduction is the field fixed by $\Ker\nu$.

\medskip
The proof of theorem~\ref{descthm} is written in~\cite{Vo} for elliptic curves but it clearly applies to abelian varieties without change. If we replace $T_p({\mathcal A})$ by $V_p({\mathcal A})$ the conclusion is that ${\mathcal A}$ is $K$-isogenous to an abelian variety defined over $F$.

\medskip
Consider now the following situation. Let $\Qp \subseteq L \subseteq K$ be finite extensions such that $K/\Qp$ is Galois, $L/\Qp$ is totally ramified, and $K/L$ is unramified. Hence 
$$\gal(K/\Qp) = I(K/\Qp)\rtimes \gal(K_0/\Qp)$$ 
with $K_0$ the maximal unramified subfield of $K$. Suppose we are given an abelian variety ${\mathcal A}$ over $L$ having good reduction. It then follows from theorem~\ref{descthm} together with~\cite{Vo} Lemma 4.6 that ${\mathcal A}$ is defined over $\Qp$ if and only if the action of $G_K$ on $T_p({\mathcal A})$ extends to an action of $G$ such that: 
\begin{itemize}
\item[(a)] the action of $G_L$ coincides with the natural action on $T_p({\mathcal A})$, and
\item[(b)] the action of $G_{K_0}$ induces an action of $I(K/\Qp)=\gal(K/K_0)$ that comes from $k$-automorphisms of the special fibre of ${\mathcal A}\times_L K$.
\end{itemize}

\subsection{Duality and polarisations}
\label{sec:dual}
We keep the hypotheses and notations of section~\ref{sec:inertia}. Let ${\mathcal A}_0^{\vee}$ be the dual abelian variety of ${\mathcal A}_0$. It also acquires good reduction over $K$ (and $L$), being $\Qp$-isogenous to ${\mathcal A}_0$, thus yielding an injective morphism
$$\nu^{\prime} : I(K/\Qp) \hookrightarrow \aut_k(A^{\vee}).$$ 
Consider the descent system of $K$-isomorphisms furnished by the Weil criterion for ${\mathcal A}_0^{\vee}$
$$\hat{f}_g : {\mathcal A}^{\vee} \xrightarrow{\sim} ({\mathcal A}^{\vee})^g  = ({\mathcal A}^g)^{\vee} , \quad g \in \gal(K/\Qp).$$
On the other hand, by dualising the system $f_g$ we had for ${\mathcal A}_0$ and taking inverses, we obtain a set of $K$-isomorphisms for ${\mathcal A}_0^{\vee}$ that satisfies the cocycle condition. By unicity the two systems must match, that is, $\hat{f}_g = \bigl( (f_g)^{\vee}\bigr)^{-1}$ for all $g \in \gal(K/\Qp)$. Hence  
$$\nu^{\prime}(g) = (\nu(g)^{\vee})^{-1} \quad \text{ for all } g\in I(K/\Qp).$$
Now take any $\Qp$-polarisation $\Lambda_0 : {\mathcal A}_0 \rightarrow {\mathcal A}_0^{\vee}$. Let $\Lambda = \Lambda_0\times_{\Qp}K : {\mathcal A} \rightarrow {\mathcal A}^{\vee}$ and $\lambda : A \rightarrow A^{\vee}$ obtained from $\Lambda$ by taking special fibres of N\'eron models.

\begin{lem}
\label{pgrpollem}
We have $\nu\bigl(I(K/\Qp)\bigr) \subseteq \aut_k(A,\lambda)$ for any polarisation $\lambda$ on $A$ deduced from one on ${\mathcal A}_0$.
\end{lem}

\begin{proof}
Since $\Lambda$ comes from the $\Qp$-polarisation $\Lambda_0$, it is compatible with the descent systems of ${\mathcal A}_0$ and ${\mathcal A}_0^{\vee}$, that is, $\Lambda^g  f_g = \hat{f}_g \Lambda$ for all $g \in \gal(K/\Qp)$. Knowing that $\hat{f}_g = \bigl( (f_g)^{\vee}\bigr)^{-1}$ we obtain $(f_g)^{\vee} \Lambda^g f_g = \Lambda$ for all $g$, which in turn gives, by restriction to the inertia subgroup and taking special fibres, $\nu(g)^{\vee} \lambda \, \nu(g) = \lambda$ for all $g \in I(K/\Qp)$.
\end{proof}

\section{Abelian varieties over finite fields}
\label{sec:abvarfinite}

\subsection{The Honda-Tate theory}
\label{sec:HT}
Let $A$ be an abelian variety over $k=\F_q$ with $q=p^s$ and let ${\rm Frob}_A$ be its Frobenius endomorphism. The algebra $\Endom^{\circ}_{k}(A)=\Q\otimes_{\Z}\Endom_{k}(A)$ is semisimple with centre $\Q({\rm Frob}_A)$, call it the Honda-Tate algebra associated to $A/k$. The characteristic polynomial $P_A=\pchar({\rm Frob}_A)$ is monic of degree $2\dim A$ with coefficients in $\Z$. The splitting of $A$ up to $k$-isogeny into powers of pairwise nonisogenous simple varieties corresponds to the Wedderburn decomposition of its Honda-Tate algebra, which itself is given by the factorisation of $P_A$ into $\Q$-irreducibles. If $A\sim B^n$ ($k$-isogenous) then $\Endom^{\circ}_{k}(A)=M_n\bigl(\Endom^{\circ}_{k}(B)\bigr)$ and if $B$ is simple then $\Endom^{\circ}_{k}(B)$ is a skewfield. Tate's theorems assert that for any $A$ and $B$ we have canonical isomorphisms 
$$\Q_{\l}\otimes_{\Z}\Hom_k(A,B)\simeq \Hom_{_{\text{Frob}}}\bigl(V_{\l}(A),V_{\l}(B)\bigr)$$
for each $\l \neq p$ (\cite{Ta1}) and, letting $D(A)=\D\bigl(A(p)\bigr)$, 
$$\Q_p\otimes_{\Z}\Hom_k(A,B)\simeq \Hom_{\varphi}\bigl(D(B),D(A)\bigr)$$
(see~\cite{Wa-Mi} for a proof). In particular $A\sim B$ if and only if $P_A=P_B$.

\medskip
A $q$-Weil number is an algebraic integer such that all its conjugates have absolute value $\sqrt{q}$ in $\C$. It is known that the roots of $P_A$ are $q$-Weil numbers. Thus the conjugacy class of a root $\pi$ of $P_A$ determines the isogeny class of $A$. Conversely, if $\pi$ is a $q$-Weil number then by Honda's theorem~\cite{Ho-Ta} there exists a simple abelian variety $A$ over $k$ such that $P_A=\pmin(\pi)^{\delta}$, where $\delta$ is the degree of the division algebra $D=\Endom^{\circ}_{k}(A)$ with centre $F=\Q(\pi)$. Furthermore the invariant of $D$ at a prime $v$ of $F$ is given by 
$$\inv(F_v\otimes_FD) = 
\begin{cases}
0 & \text{if $v\mid \l$ with $\l\neq p$ a rational prime} \\ 
\frac{1}{2} & \text{if $v$ is real} \\
\frac{f_v\ord_v(\pi)}{s} & \text{if $v\mid p$, where $f_v=f(F_v/\Qp)$ is the residue degree.}
\end{cases}$$
In particular, for $v$ lying above $p$, $F_v\otimes_FD$ is Brauer equivalent to the cyclic algebra 
$$\Bigl( K_0F_v/F_v , \sigma^{f_v}, \pi^{f_v/(s,f_v)} \Bigr) = \Bigl( K_0F_v/F_v, \sigma^{(s,f_v)}, \pi \Bigr)$$
with $K_0=\Frac W(\F_{p^s})$ and $(s,f_v)=\gcd(s,f_v)$. If $k=\F_p$ and $\pi$ has no real conjugate then $\Endom^{\circ}_{k}(A)=F$. The characteristic polynomial of a simple $A/\F_{p^s}$ has real roots in only two cases:
\begin{itemize}
\item[(a)] $s$ is even and $P_A(X)=(X+ p^{s/2})^2$ or $(X- p^{s/2})^2$ (corresponding to a quadratic twist), then $A$ is a supersingular elliptic curve with $\Endom^{\circ}_{\F_{p^s}}(A)=D_{p,\infty}=$ the quaternion algebra over $\Q$ ramified only at $p$ and $\infty$ 
\item[(b)] $s$ is odd and $P_A(X)=(X^2-p^s)^2$, then $A$ is an abelian surface with $\Endom^{\circ}_{\F_{p^s}}(A)={\mathcal D}_{\infty}=$ the quaternion algebra over $\Q(\sqrt{p})$ ramified only at the two infinite primes. Over a quadratic extension $A$ becomes isogenous to the product of two supersingular elliptic curves of the second type above. 
\end{itemize}
So the invariant $\delta$ attached to the simple abelian variety $A/k$, or equivalently to the $q$-Weil number $\pi$ such that $P_A(\pi)=0$, is given by
$$\delta = 
\begin{cases}
2 & \text{if $\pmin(\pi)$ has a real root} \\
\lcm \biggl[ \ord_{\Q/\Z} \Bigl( \frac{f_v\ord_v(\pi)}{s} \Bigr) \ ,\  v\mid p \text{ in } F \biggr] & \text{otherwise.}
\end{cases}$$

\begin{defn}
\label{pWpoldefn}
A monic polynomial $P\in \Z[X]$ is a {\em $p$-Weil polynomial} if all its roots in $\Qa$ are $p$-Weil numbers and its valuation at $X^2-p$ is even.
\end{defn}

Equivalently, all the roots are $p$-Weil numbers and $P(0)=p^{\deg P/2}$. Thus 
$$P(X)=(X^2-p)^{2n_0}P_1(X)^{n_1}\dots P_m(X)^{n_m}$$ 
with $P_i$ irreducibles in $\Z[X]$ having no real root. The characteristic polynomial of an abelian variety over $\F_p$ is a $p$-Weil polynomial, and for each $p$-Weil polynomial $P$ there is an $A/\F_p$ such that $P_A=P$. Specifically, there is a simple abelian surface $A_0$ with $P_{A_0}(X)=(X^2-p)^2$ and for each $i$ there is a simple $A_i$ with $P_{A_i}=P_i$, so that we can take 
$$A = \prod_{0\leq i\leq m} A_i^{n_i}.$$
We thus get a convenient correspondance between $p$-Weil polynomials and isogeny classes of abelian varieties over $\F_p$.

\medskip
We conclude this section by an elementary remark on the behaviour under base change. For $P\in \Q[X]$ splitting over $\Qa$ as $P(X)=\prod_i (X-\lambda_i)$ and $s$ a positive integer, put
$$P^{(s)}(X) \underset{\text{def}}{=} \; \prod_i (X-\lambda_i^s).$$
One checks that $P^{(s)} \in \Q[X]$, and if $P$ is irreducible then $P^{(s)}$ is a power of an irreducible in $\Q[X]$. Now for $A=A_0\times_{\F_p}\F_{p^s}$ we have $P_A=P_{A_0}^{(s)}$, consequently if $A_0/\F_p$ is simple then $A\sim B^n$ with $B/\F_{p^s}$ simple.

\subsection{Honda-Tate algebras}
\label{sec:HTalg}   
In the sequel we will need to deal with some explicit aspects of the Honda-Tate theory, which is the purpose of this section. Lemma~\ref{splitlem} furnishes a cyclotomic extension unramified at $p$ that splits a given simple Honda-Tate algebra. We then investigate under which condition this algebra contains a certain cyclic subalgebra that naturally appears in section~\ref{sec:tamegalpairs}. We state this condition in terms of some arithmetic invariant in proposition~\ref{embedprop}.

\medskip
We begin with a result on subalgebras, known to the experts as Schofield's Lemma, that will be used several times in this paper without comment. The proof we present here for the reader's convenience is due to E. Frossard. Let $K$ be any field and $\mathsf{A}$, $\mathsf{B}$ two central simple algebras over $K$. Denote by $[\mathsf{A}]$ the class of $\mathsf{A}$ in the Brauer group of $K$. 

\begin{lem}[\cite{Sc} 5, Lemma 9.1]
\label{sbalglem}
There is a $K$-algebra embedding  $\mathsf{A} \hookrightarrow \mathsf{B}$ if and only if $\ind\bigl([\mathsf{B}]-[\mathsf{A}]\bigr)\deg\mathsf{A}$ divides $\deg\mathsf{B}$.  
\end{lem}

\begin{proof}
Assume $\mathsf{A}\hookrightarrow \mathsf{B}$ and identify $\mathsf{A}$ with its image. Obviously $\deg\mathsf{A}$ divides $\deg\mathsf{B}$. Let $\mathsf{C}$ be the centraliser of $\mathsf{A}$ in $\mathsf{B}$. Then $\mathsf{A}\otimes_K\mathsf{C}=\mathsf{B}$, so $[\mathsf{C}]= [\mathsf{B}]-[\mathsf{A}]$ and $\ind\bigl([\mathsf{C}]\bigr)$ divides $\deg\mathsf{C} = \deg\mathsf{B}/\deg\mathsf{A}$. 

Conversely, assume $\deg\mathsf{B}=n \ind\bigl([\mathsf{B}]-[\mathsf{A}]\bigr)\deg\mathsf{A}$ for some integer $n$. Let $\mathsf{D}$ be the skewfield such that $[\mathsf{D}]= [\mathsf{B}]-[\mathsf{A}]$ and put $\mathsf{B}^{\prime}=\mathsf{A}\otimes_K M_n(\mathsf{D})$. Then $[\mathsf{B}^{\prime}]=[\mathsf{A}]+[\mathsf{D}]=[\mathsf{B}]$ and $\deg\mathsf{B}^{\prime}=n\deg\mathsf{D}\deg\mathsf{A}=n\ind\bigl([\mathsf{B}]-[\mathsf{A}]\bigr)\deg\mathsf{A}=\deg\mathsf{B}$, hence $\mathsf{B}\simeq \mathsf{B}^{\prime}$.
\end{proof}

\begin{cor}
\label{locglobcor}
Let $K$ be a number field. Then $\mathsf{A} \hookrightarrow \mathsf{B}$ if and only if $K_v\otimes_K\mathsf{A} \hookrightarrow K_v\otimes_K\mathsf{B}$ for all places $v$ of $K$.
\end{cor}

\begin{proof}
Since $K$ is a number field the global index is the least common multiple of the local indices, whereas degrees are unchanged. 
\end{proof}

\medskip
Let $B/\F_{p^s}$ be a simple abelian variety with $p^s$-Weil number $\pi$, and $D$ its Honda-Tate division algebra with centre $F=\Q(\pi)$ and index $\delta$. Let $\zeta_r\in \Qa$ be a primitive $r$th root of unity.

\begin{lem}
\label{splitlem}
Assume $r\geq3$ is an integer prime to $p$ such that $s$ is the order of $p$ in $(\Z/r\Z)^{\times}$. Then $F(\zeta_r)$ splits $D$. Furthermore $\delta$ divides $[F(\zeta_r):F]$ and $n=[F(\zeta_r):F]/\delta$ is the smallest integer such that $F(\zeta_r)$ embeds in $M_n(D)$.
\end{lem}

\begin{proof}
Let $v$ be a prime of $F$ lying above $p$ and $f_v$ the residue degree of $F_v$. Our assumption on $r$ implies that $K_0F_v=F(\zeta_r)_v$ with $K_0=\Frac W(\F_{p^s})$. We have $[K_0F_v:F_v]=s/(s,f_v)$, where $(s,f_v)=\gcd(s,f_v)$, thus 
$$[K_0F_v:F_v]\, \inv(F_v\otimes_FD) = \frac{f_v}{(s,f_v)}\ord_v(\pi) \in \Z.$$
So $F(\zeta_r)_v$ splits $F_v\otimes_FD$. At primes above $\l\neq p$ we are in a split situation to start with. At an archimedian prime $\infty$ we have $F(\zeta_r)_{\infty}=\C$ since $r\geq3$. Therefore $F(\zeta_r)$ splits $D$. The rest follows from the index reduction theorem, see~\cite{Al} IV, $\S\S$9,10. Then $F(\zeta_r)$ is a maximal subfield of $M_{[F(\zeta_r):F]/\delta}(D)$.
\end{proof}

\medskip
Fix $r\geq3$ prime to $p$ such that $s$ is the order of $p$ in $(\Z/r\Z)^{\times}$. Let $m$ be a multiple of $\delta$ and consider the Honda-Tate algebra 
$$M_{m/\delta}(D) = \Endom_{\F_{p^s}}^{\circ}(A) \quad \text{ with } A \sim B^{m/\delta}.$$
We assume that $\gal\bigl(F(\zeta_r)/F\bigr)$ contains an automorphism $\sigma_p$ acting as $\zeta_r \mapsto \zeta_r^p$ (i.e. $F\cap \Q(\zeta_r)$ is contained in the subfield of $\Q(\zeta_r)$ fixed by $\zeta_r \mapsto \zeta_r^p$). Let $L$ be the field fixed by $\sigma_p$. We want to state necessary and sufficient conditions for $M_{m/\delta}(D)$ to contain a subalgebra isomorphic to the cyclic algebra of degree $s$
$$C(\pi) \underset{\text{def}}{=} \bigl(F(\zeta_r)/L,\sigma_p,\pi\bigr).$$
Lemma~\ref{splitlem} shows that $F(\zeta_r)$ embeds in $M_{m/\delta}(D)$ if and only if $[F(\zeta_r):F]$ divides $m$. The algebra $M_{[F(\zeta_r):F]/\delta}(D)$ is then a crossed product $\bigl( F(\zeta_r)/F,\mathsf{c} \bigr)$ for some $2$-cocycle $\mathsf{c}\in H^2\bigl(\gal\bigl(F(\zeta_r)/F\bigr),F(\zeta_r)^{\times}\bigr)$. Pick an $a\in L^{\times}$ such that
$${\rm Res}_{\langle \sigma_p\rangle}\bigl( F(\zeta_r)/F,\mathsf{c} \bigr) \simeq \bigl(F(\zeta_r)/L,\sigma_p,a\bigr) \underset{\text{def}}{=}C(a)$$
(unique up to a norm from $F(\zeta_r)$ to $L$); this cyclic algebra is Brauer equivalent to $L\otimes_FD$, being the centraliser of $L$ in $\bigl( F(\zeta_r)/F,\mathsf{c} \bigr)$. Note that $C(\pi)$ embeds in $M_{m/\delta}(D)$ if and only if it embeds in $M_{m/[F(\zeta_r):F]}\bigl(C(a)\bigr)$, the centraliser of $L$ in $M_{m/\delta}(D)$. 

\medskip
For each place $v$ of $L$ define 
$$n_v(r;\pi) \underset{\text{def}}{=} 
\begin{cases}
\ord \Bigl( \pi \in  L_v^{\times} {\big /} N_{F(\zeta_r)_v/L_v}\bigl(F(\zeta_r)_v^{\times}\bigr) \Bigr) & \text{if $v\mid \l \neq p$ and $\l\mid r$} \\ 
2 & \text{if $L_v=\R$ and $\pi= p^{s/2}$ with $s$ even} \\
1 & \text{otherwise.}
\end{cases}$$
This integer is actually the order of $\pi/a$ in $L_v^{\times} {\big /} N_{F(\zeta_r)_v/L_v}\bigl(F(\zeta_r)_v^{\times}\bigr)$. Indeed, it follows from Honda-Tate that $\pi /a$ is a local norm at each $v\mid p$ and that $a$ is a local norm at each $v\mid \l \neq p$. Since $\pi$ is an $\l$-adic unit for $\l \neq p$, it is a local norm at $v\mid \l$ whenever $F(\zeta_r)/L$ is unramified at $v$, which is the case if $\l$ does not divide $r$. Finally, $L$ has a real place if and only if $s$ is even, $\pi= \pm p^{s/2}$, and $p^{s/2}\equiv -1 \bmod r$; then $F=\Q$, $D=D_{p,\infty}$ and $a$ is not a norm from $\C$ to $\R$. Now put 
$$n(r;\pi) \underset{\text{def}}{=} \lcm \Bigl(n_v(r;\pi), \text{ $v$ place of $L$} \Bigr).$$
The extension $F(\zeta_r)/L$ being cyclic $n(r;\pi)$ is the order of $\pi/a$ in $L^{\times} {\big /} N_{F(\zeta_r)/L}\bigl(F(\zeta_r)^{\times}\bigr)$.

\begin{prop}
\label{embedprop}
There is an embedding $\bigl(F(\zeta_r)/L,\sigma_p,\pi\bigr) \hookrightarrow M_{m/\delta}(D)= \Endom_{\F_{p^s}}^{\circ}(A)$ if and only if $[F(\zeta_r):F]n(r;\pi)$ divides $m$.
\end{prop}

\begin{proof}
By lemma~\ref{splitlem} $[F(\zeta_r):F]$ must divide $m$. According to corollary~\ref{locglobcor} we have to show that, for each place $v$ of $L$,   
$$(\iota_v) \quad \quad \quad  L_v\otimes_L C(\pi) \hookrightarrow L_v\otimes_L M_{m/[F(\zeta_r):F]}\bigl(C(a)\bigr) $$
holds if and only if $n_v(r;\pi)$ divides $m/[F(\zeta_r):F]$. For $v\mid p$ we have $F_v=L_v$ and $L_v\otimes_L C(\pi) \simeq \bigl( K_0F_v/F_v,\sigma^{(s,f_v)},\pi \bigr) \simeq L_v\otimes_L C(a)$, so $(\iota_v)$ holds. For $v\mid \l \neq p$ we have
$$L_v\otimes_L M_{m/[F(\zeta_r):F]}\bigl(C(a)\bigr) \simeq M_{m/[L:F]}(L_v) \quad \text{whereas} \quad L_v\otimes_L C(\pi) \simeq M_{s/n_v(r;\pi)}(D_v)$$
with $D_v$ the division algebra of degree $n_v(r;\pi)$ over $L_v$. Hence $(\iota_v)$ holds if and only if $n_v(r;\pi)$ divides $m/[F(\zeta_r):F]$. Finally, let $v$ be archimedean. If $L_v=\C$ then $(\iota_v)$ holds. If $L_v=\R$ then $s$ is even, $\pi= \pm p^{s/2}$, $D=D_{p,\infty}$, and $\R\otimes_L C(a) \simeq M_{s/2}(\H)$. For $\pi= - p^{s/2}$ we have $\R\otimes_L C(\pi) \simeq M_{s/2}(\H)$ and $(\iota_v)$ holds. For $\pi= p^{s/2}$ we have $\R\otimes_L C(\pi) \simeq M_{s}(\R)$, so $(\iota_v)$ holds if and only if $n_v(r;\pi)=2$ divides $m/[F(\zeta_r):F]$.
\end{proof}

\section{Galois pairs}
\label{sec:galpairs}
Unless otherwise specified the following notations will be in force throughout the rest of the paper. For an abelian variety $A_0/\F_p$ let $\mathsf{f}_0={\rm Frob}_{A_0}$, $A=A_0\times_{\F_p}\F_{p^s}$, $\mathsf{f}={\rm Frob}_{A}= \mathsf{f}_0^s$, $K_0=\Frac W(\F_{p^s})$, $D(A)=\D\bigl(A(p)\bigr)$ and $\Delta(A)=\W\bigl(D(A)\bigr)$.

\medskip
Recall that $D(A)$ is a $K_0$-vector space of dimension $2\dim A$ together with a $\sigma$-semilinear bijection $\varphi: D(A)\xrightarrow{\sim}D(A)$ such that $\pchar(\varphi^s)=\pchar(\mathsf{f})$, and $D(A_0)=\D\bigl(A_0(p)\bigr)$ is a $\Qp$-vector space of the same dimension together with a $\Qp$-linear isomorphism $\varphi_0: D(A_0)\xrightarrow{\sim}D(A_0)$ such that $\pchar(\varphi_0)=\pchar(\mathsf{f}_0)$. We actually have 
$$D(A)=K_0\otimes_{\Qp}D(A_0) \quad \text{and} \quad \varphi = \sigma \otimes \varphi_0.$$ 
Thus $D(A)$ has a natural structure of a $\bigl(\varphi,\gal(K_0/\Qp)\bigr)$-module. Indeed, $\sigma(\lambda \otimes x_0)= \sigma(\lambda)\otimes x_0$ for $\lambda \in K_0$, $x_0\in D(A_0)$ defines a $\sigma$-semilinear action of $\gal(K_0/\Qp)$ commuting with $\varphi$, and then $D(A_0)=\{ x\in D(A) \mid \sigma x=x\}$. We have 
$$\Qp \otimes_{\Q} \Endom_{\F_p}^{\circ}(A_0) \simeq \Endom_{\varphi_0}\bigl(D(A_0)\bigr) \simeq \Endom_{\varphi,\gal(K_0/\Qp)}\bigl(D(A)\bigr).$$
We want to mimic the situation considered in section~\ref{sec:inertia} where a $\bigl(\varphi,\gal(K/\Qp)\bigr)$-module structure was obtained on $D(A)$ in a context of potential good reduction, with $K$ a suitable Galois extension of $\Qp$ containing $K_0$. We first define in~\ref{sec:defn} the appropriate notion, namely Galois pairs. In~\ref{sec:repnsgalpairs} we check that these objects yield the desired representations, be it $\bigl(\varphi,\gal(K/\Qp)\bigr)$-modules or Weil representations. We then specialise to the tame case in~\ref{sec:tamegalpairs} where we prove a decomposition result (proposition~\ref{decompprop}).

\subsection{Definition of Galois pairs}
\label{sec:defn}
We start with an abelian variety $A$ over $\F_{p^s}$ and consider the set of $\F_{p^s}$-isomorphisms from $A$ to its various $\gal(\F_{p^s}/\F_p)$-twists 
$$\ist_{\F_{p^s}}(A) \underset{\text{def}}{=} \{\, \F_{p^s}\text{-isomorphisms }A\xrightarrow{\sim}A^{g},\; g\in \gal(\F_{p^s}/\F_p) \,\}.$$
For $\psi_g : A\xrightarrow{\sim}A^{g}$ and $\psi_h : A\xrightarrow{\sim}A^{h}$ with $g,h\in \gal(\F_{p^s}/\F_p)$, set 
$$\psi_g * \psi_h  \underset{\text{def}}{=} (\psi_g)^{h}\circ \psi_h \; : A\xrightarrow{\sim}A^{gh}=A^{hg}.$$
This defines a group structure on $\ist_{\F_{p^s}}(A)$ with identity $\id_A$ (the inverse of $\psi_g$ is $(\psi_g^{-1})^{g^{-1}}$). $\aut_{\F_{p^s}}(A)$ is an invariant subgroup, the restriction of $*$ being the usual group law, and the quotient identifies with a subgroup of $\gal(\F_{p^s}/\F_p)$. Write $\psi^{[n]}=\psi * \dots *\psi \text{ ($n$ times)}$ for $n\in \Z$.

Now assume $A=A_0\times_{\F_p}\F_{p^s}$ for an abelian variety $A_0/\F_p$. Then the Weil criterion furnishes a distinguished element $f_{\sigma}\in \ist_{\F_{p^s}}(A)$, namely 
$$f_{\sigma}=\id_{A_0}\otimes \sigma : A\xrightarrow{\sim}A^{\sigma}$$ 
generating a subgroup of order $s$. The map $f_{\sigma}^{[n]}\mapsto \sigma^n$ induces an isomorphism $\langle f_{\sigma}\rangle \simeq \gal(\F_{p^s}/\F_p)$ which provides a section for the exact sequence
$$1 \longrightarrow \aut_{\F_{p^s}}(A) \longrightarrow \ist_{\F_{p^s}}(A) \longrightarrow \gal(\F_{p^s}/\F_p) \longrightarrow 1.$$
Therefore in this situation $\ist_{\F_{p^s}}(A)= \aut_{\F_{p^s}}(A) \rtimes \langle f_{\sigma}\rangle$.

\begin{defn}
\label{galabdefn} 
Let $K/\Qp$ be finite Galois with residue field $\F_{p^s}$. Let $A_0/\F_p$ be an abelian variety and $\Gamma$ a finite subgroup of $\aut_{\F_{p^s}}(A)$. Say $(A_0,\Gamma)$ is a {\em Galois pair} for $K/\Qp$ if there exists an injective (anti-)morphism 
$$\nu \; : \; 
\begin{cases}
\gal(K/\Qp) & \hookrightarrow  \quad \ist_{\F_{p^s}}(A)   \\
g           & \mapsto  \quad \nu(g) : A\xrightarrow{\sim}A^{\bar{g}} \quad \text{ with }\bar{g}=g\bmod I(K/\Qp)
\end{cases}$$
satisfying the following conditions:
\begin{itemize}
\item[(i)] $\nu \bigl(I(K/\Qp)\bigr)= \Gamma$, and
\item[(ii)]$\nu(\omega)=f_{\sigma}$ for some $\omega\in\gal(K/\Qp)$ such that $\bar{\omega}=\sigma$.                  
\end{itemize}
\end{defn}

We write $(A_0,\Gamma)$ instead of $(A_0,\Gamma,\nu)$ in order to ease the terminology and notations. Note that we always have $\nu \bigl(I(K/\Qp)\bigr) \subseteq \aut_{\F_{p^s}}(A)$; condition $(i)$ imposes the given $\Gamma$ to be the image under $\nu$ of the inertia subgroup. Condition $(ii)$ is equivalent to the commutativity of the diagram
$$\xymatrix{
1 \ar[r]  &  I(K/\Qp) \ar[r] \ar@{^{(}->}[d]_{\nu}  &  \gal(K/\Qp) \ar[r] \ar@{^{(}->}[d]_{\nu}   &  \gal(\F_{p^s}/\F_p) \ar[r] \ar@{=}[d]  &   1  \\
1 \ar[r]  &  \aut_{\F_{p^s}}(A)   \ar[r]            &  \ist_{\F_{p^s}}(A)  \ar[r]                 &   \gal(\F_{p^s}/\F_p) \ar[r]            &   1  }$$
In particular we see that $\gal(K/\Qp)= I(K/\Qp)\rtimes \gal(K_0/\Qp)$, hence $K$ is an unramified extension of a totally ramifed extension $L$ of $\Qp$ (compare with section~\ref{sec:inertia}). Observe that we may always assume $s$ to be minimal with respect to $\Gamma$, that is, that $\F_{p^s}$ is the smallest field of definition for $\Gamma$; then $K$ is the Galois closure of $L$.

\begin{ex}
Let $e\in\{ 3,4,6\}$ and assume $(p,e)=1$. Let $E_0/\F_p$ be an ordinary elliptic curve with $\zeta_e\in \aut_{\F_p}(E_0)$. Then $p\equiv 1\bmod e$ and $(E_0,\langle\zeta_e\rangle)$ is a Galois pair ($s=1$) for a totally ramified extension $K/\Qp$ of degree $e$. If $E_0/\F_p$ is supersingular with $\zeta_e\in \aut_{\Fpa}(E_0)$ then $p\equiv -1\bmod e$, and $(E_0,\langle\zeta_e\rangle)$ is a Galois pair ($s=2$) for $K/\Qp$ with $K$ a totally ramified extension of $\Frac W(\F_{p^2})$ of degree $e$. 
\end{ex}

Given $A_0$, not all pairs $(A_0,\Gamma)$ with $\Gamma$ a finite subgroup of $\aut_{\Fpa}(A_0)$ are Galois. 

\begin{ex}
Let $A_0={B_0}^n$ so that $\aut_{\F_p}(A_0)$ contains a subgroup isomorphic to the permutation group ${\mathcal S}_n$. Let $\Gamma$ be a cyclic subgroup of ${\mathcal S}_n$ of order $e$ prime to $p$. Then $(A_0,\Gamma)$ is a Galois pair if and only if $p\equiv 1\bmod e$.
\end{ex}

\medskip
A morphism of Galois pairs $\psi_0 : (A_0,\Gamma,\nu) \rightarrow (A_0^{\prime},\Gamma^{\prime},\nu^{\prime})$ is an $\F_p$-morphism of abelian varieties $\psi_0 : A_0\rightarrow A_0^{\prime}$ such that the diagram 
$$\xymatrix{
A \ar[r]^{\psi} \ar[d]_{\nu(g)}^{\wr}  &  A^{\prime} \ar[d]^{\nu^{\prime}(g)}_{\wr}   \\
A^{\bar{g}} \ar[r]_{\psi^{\bar{g}}}    &  {A^{\prime}}^{\bar{g}}                }$$
commutes for all $g \in \gal(K/\Qp)$. There is an obvious notion of isogeny of Galois pairs, just require in addition $\psi_0$ to be an isogeny.

\subsection{The representations associated to Galois pairs}
\label{sec:repnsgalpairs}
The definition of a Galois pair $(A_0,\Gamma)$ for $K/\Qp$ has been taylored so to provide a $\bigl(\varphi,\gal(K/\Qp)\bigr)$-module structure on $D(A)$, via $\nu$. Indeed, passing to $p$-divisible groups and applying the functor $\D$ we get 
$$\xymatrix{
\gal(K/\Qp) \ar@{^{(}->}[r]^{\nu} \ar@{^{(}->}[dr]_{\nu}  &  \ist_{\F_{p^s}}(A) \ar@{^{(}->}[d]^{\text{functor $D$}} \\
                                                          &  \ist_{\F_{p^s}}^*\bigl(D(A)\bigr)                       }$$
where $\ist_{\F_{p^s}}^*\bigl(D(A)\bigr)$ is the group of $\varphi$-modules isomorphisms $D(A^g)\xrightarrow{\sim}D(A)$ for all $g \in \gal(\F_{p^s}/\F_p)$, the group law being defined in a similar way as the one for $\ist_{\F_{p^s}}(A)$. Note that this time the diagonal $\nu$ is a genuine group morphism, thanks to contravariancy. This obviously defines a faithful $\sigma$-semilinear action of $\gal(K/\Qp)$ on $D(A)$. Since $f_{\sigma}*\mathsf{f}_0 = f_{\sigma}\circ \mathsf{f}_0$ commutes with all elements in $\ist_{\F_{p^s}}(A)$ and
$$D(f_{\sigma}*\mathsf{f}_0) = \sigma \otimes \varphi_0 = \varphi$$
this action commutes with $\varphi$. We thus obtain the $\bigl(\varphi,\gal(K/\Qp)\bigr)$-module structure on $D(A)$ associated to the Galois pair $(A_0,\Gamma)$. Obviously two isogenous Galois pairs for $K/\Qp$ yield isomorphic $\bigl(\varphi,\gal(K/\Qp)\bigr)$-modules. 

\begin{rem}
\label{Tatewithdescentrmk}
Let $\Endom_{\F_p,\Gamma}^{\circ}(A_0)$ be the subalgebra of $\Endom_{\F_p}^{\circ}(A_0)$ consisting of elements commuting with all elements in $\Gamma$ (after base change to $\F_{p^s}$). We have 
$$\Qp \otimes_{\Q} \Endom_{\F_p,\Gamma}^{\circ}(A_0) \simeq \Endom_{\varphi,\gal(K/\Qp)}\bigl(D(A)\bigr).$$
\end{rem}

\medskip
We now construct the associated Weil representation. Recall that $\upsilon : W \twoheadrightarrow \Z$ is the morphism sending a lifting $\sigma_0$ of the absolute Frobenius to $1$. Put
$$\rho \; : \; 
\begin{cases}
W & \rightarrow  \quad   \Endom^{\circ}_{\F_{p^s}}(A)^{\times} \\
w & \mapsto      \quad   \Bigl( (f_{\sigma}*\mathsf{f}_0)^{[\upsilon(w)]} \Bigr)^{-1} \nu(w\bmod W_K)
\end{cases}$$
One checks that $\rho$ is a well-defined antimorphism satisfying $\rho(I)=\Gamma$ and $\rho(\sigma_0^{-1})=\mathsf{f}_0$. If we let $M=K\Qp^{\text{un}}$ and ${\rm Weil}(M/\Qp)= W/I_K =I(K/\Qp)\rtimes\langle\sigma_0\rangle$ then $\rho$ yields an injection
$$\rho : {\rm Weil}(M/\Qp) \hookrightarrow \aut_{\F_{p^s}}(A) \rtimes \langle \mathsf{f}_0 \rangle \subseteq \Endom^{\circ}_{\F_{p^s}}(A)^{\times}$$
with image $\Gamma \rtimes \langle \mathsf{f}_0 \rangle$. We thus obtain the associated Weil representation 
$$\xymatrix{
{\rm Weil}(M/\Qp) \ar@{^{(}->}[r]^{\rho} \ar@{^{(}->}[dr]_{\rho}  &  \Endom^{\circ}_{\F_{p^s}}(A)^{\times} \ar@{^{(}->}[d]^{\text{functor $\Delta$}} \\
                                                                  &  \aut_{K_0}\bigl(\Delta(A)\bigr)                                                 }$$

\medskip
So this is how we get $\rho$ from $\nu$. Conversely, suppose we are given a $\rho$ as above. We define an injective antimorphism by
$$\nu \; : \; 
\begin{cases}
\gal(K/\Qp)   & \hookrightarrow \quad  \ist_{\F_{p^s}}(A) = \aut_{\F_{p^s}}(A)\rtimes \langle f_{\sigma}\rangle \\
 g            & \mapsto         \quad  \eta \bigl( \rho(\tilde{g}) \bmod \langle \mathsf{f} \rangle \bigr)
\end{cases}$$
where $\tilde{g}$ is any lifting of $g$ in ${\rm Weil}(M/\Qp)$ and $\eta$ is $(\mathsf{f}_0\bmod \langle \mathsf{f} \rangle)^n\gamma \mapsto f_{\sigma}^{[n]}*\gamma$ for $\gamma\in\aut_{\F_{p^s}}(A)$. We then have a commutative diagram
$$\xymatrix{
{\rm Weil}(M/\Qp) \ar@{->>}[dd]_{\bmod W_K/I_K}  \ar@{^{(}->}[r]^{\rho}  &  \aut_{\F_{p^s}}(A) \rtimes \langle \mathsf{f}_0 \rangle \ar@{->>}[d]^{\bmod \langle \mathsf{f} \rangle}  \\
    &  \aut_{\F_{p^s}}(A) \rtimes \langle \mathsf{f}_0 \bmod \langle \mathsf{f} \rangle \rangle  \ar[d]^{\eta}_{\wr}    \\
\gal(K/\Qp) \ar@{^{(}->}[r]^{\nu}  &  \ist_{\F_{p^s}}(A) }$$
Hence these constructions are inverse one to the other, and, when passing to representations, $\nu \mapsto \rho$ corresponds to $D \mapsto \W(D)=\Delta$ whereas $\rho \mapsto \nu$ corresponds to a ``section'' $\Delta \mapsto D$.

\begin{rem}
\label{sminrem}
The following are equivalent:
\begin{itemize}
\item[(i)] $\F_{p^s}$ is the smallest field of definition for $\Gamma$
\item[(ii)] $K/\Qp$ is a Galois extension of minimal degree with $I(K/\Qp)\simeq \Gamma$
\item[(iii)] $K/\Qp$ is a Galois extension of minimal degree over which $\Delta(A)$ acquires good reduction.
\end{itemize}
\end{rem}

\begin{rem}
\label{algrepnrem}
View the embedding $\Endom^{\circ}_{\F_{p^s}}(A) \hookrightarrow \Endom_{K_0}\bigl(\Delta(A)\bigr)$ as a faithful algebra representation. The action of a simple subalgebra $E$ of $\Endom^{\circ}_{\F_{p^s}}(A)$ is $2d/f\delta$ times the direct sum of the reduced representation of $E$ (tensored by $K_0$), where $d=\dim A$, $f=[Z(E):\Q]$ and $\delta=\deg E$ (\cite{Mi} Prop.12.12, replacing $V_{\l}(A)$ by $\Delta(A)$). Recall that the reduced representation of $E$ is the direct sum of the $f$ nonisomorphic representations of $E$ over $\Qa$ furnished by the $f$ distinct embeddings of $Z(E)$ in $\Qa$. It is defined over $\Q$.
\end{rem}

\subsection{Tame Galois pairs}
\label{sec:tamegalpairs}
Let $(A_0,\Gamma)$ be a Galois pair for $K/\Qp$ such that $\Gamma \subseteq \aut_{\F_{p^s}}(A)$ is of order $e$ prime to $p$. Then the extension $K/\Qp$ is tame and we say that $(A_0,\Gamma)$ is a tame Galois pair. In this case $\Gamma$ is a cyclic group generated by an element $\tau$ satisfying the relation 
$$\mathsf{f}_0\tau=\tau^p\mathsf{f}_0.$$ 
We assume $s$ to be the order of $p$ in $(\Z/e\Z)^{\times}$ since it is the smallest power of $\mathsf{f}_0$ commuting with $\tau$, so that $\F_{p^s}$ is the smallest field of definition for $\Gamma$. Consider the semisimple subalgebra of $\Endom^{\circ}_{\F_{p^s}}(A)$ generated over $\Q$ by $\mathsf{f}_0$ and $\tau$
$$\Q[\mathsf{f}_0,\tau] \subseteq  \Endom^{\circ}_{\F_{p^s}}(A).$$

\begin{defn}
\label{cyclicdefn}
The tame Galois pair $(A_0,\langle\tau\rangle)$ is {\em cyclic} if $\Q[\mathsf{f}_0,\tau]$ is a cyclic algebra.
\end{defn}

Note that $\Endom^{\circ}_{\F_{p^s}}(A)$ is central simple when $(A_0,\langle\tau\rangle)$ is cyclic, since $\pmin(\mathsf{f})$ is then irreducible in $\Q[X]$.

\begin{prop}
\label{decompprop}
A tame Galois pair is isogenous to a product of cyclic ones. A cyclic tame Galois pair $(A_0,\langle \tau \rangle)$ with $r=\ord(\tau)$ and $s=\ord \bigl(p\in (\Z/r\Z)^{\times}\bigr)$ satisfies 
$$\Q[\mathsf{f}_0,\tau] \simeq \bigl( F(\zeta_r)/L,\sigma_p,\pi \bigr)$$ 
where $\pi$ is the associated $p^s$-Weil number, $\sigma_p$ fixes $F=\Q(\pi)$, and $\sigma_p(\zeta_r)=\zeta_r^p$. 
\end{prop}

\begin{proof}
Let $(A_0,\langle \tau \rangle)$ be a tame Galois pair with $e=\ord(\tau)$ and $s=\ord \bigl(p\in (\Z/e\Z)^{\times}\bigr)$. The minimal polynomial in $\Q[X]$ annihilating $\tau$ writes as
$$\pmin(\tau) = \prod_{r} \Phi_r \quad \text{ with } r\mid e \text{ for all }r \text{ and } \lcm(\text{all }r) = e$$ 
where $\Phi_r$ is the $r$th cyclotomic polynomial. For such an $r$ the connected component of $\Ker\bigr(\Phi_r(\tau)\bigl)$ containing the identity is an abelian subvariety $A_r$ of $A$. Thus $A \sim \prod_{r} A_r$ over $\F_{p^s}$. The $A_r$'s are stable by $\tau$ which restricts to an automorphism with $\pmin$ equal to $\Phi_r$. Since $r$ is prime to $p$ we have $\Phi_r^{(p)}=\Phi_r$, so $\pmin(\tau^p) =\pmin(\tau)^{(p)}=\pmin(\tau)$. The relation $\mathsf{f}_0\Phi_r(\tau)=\Phi_r(\tau^p)\mathsf{f}_0$ implies that $\Phi_r(\tau^p)\bigl(\mathsf{f}_0(A_r)\bigr)=0$, which shows that $\mathsf{f}_0(A_r)\subseteq A_r$ and therefore each $A_r$ is defined over $\F_p$.

So assume $\pmin(\tau)=\Phi_r$ and $s=\ord \bigl(p\in (\Z/r\Z)^{\times}\bigr)$. The centre $Z(\Q[\mathsf{f}_0,\tau])$ is a product of number fields containing $\Q[\mathsf{f}]$. As above this decomposition splits $A$ into the product of abelian varieties over $\F_{p^s}$ each of which is stable by $\tau$ and by $\mathsf{f}_0$, hence defined over $\F_p$.

Now $(A_0,\langle \tau \rangle)$ is a tame Galois pair such that 
$$\pmin(\tau)=\Phi_r \quad \text{ and } \quad Z(\Q[\mathsf{f}_0,\tau])= L \text{ is a field.}$$ 
Furthermore the \'etale algebra $\Q[\mathsf{f},\tau]$ is a field. Indeed, $\Q[\mathsf{f},\tau]= \prod_j K_j$ is a product of number fields on which $\mathsf{f}_0$ acts by conjugation $c_0 : x\mapsto \mathsf{f}_0 x \mathsf{f}_0^{-1}$. Its effect being given by $\mathsf{f} \mapsto \mathsf{f}$ and $\tau \mapsto \tau^p$  this action preserves the direct factors. Thus $\Q[\mathsf{f},\tau]^{\langle c_0\rangle}=\prod_j L_j$ with $L_j$ the subfield of $K_j$ fixed by $c_0$. But $\Q[\mathsf{f},\tau]^{\langle c_0\rangle}\subseteq L$ so $\Q[\mathsf{f},\tau]$ must be a field.

Let $\pi$ be a root of the $\Q$-irreducible $\pmin(\mathsf{f})$ and $F=\Q(\pi)$. We have an isomorphism $\Q[\mathsf{f},\tau] \simeq F(\zeta_r)$ given by $\mathsf{f}\mapsto \pi$ and $\tau\mapsto \zeta_r$, with $\zeta_r\in \Qa$ some primitive $r$th root of unity. The conjugation $c_0$ on $\Q[\mathsf{f},\tau]$ translates to an element $\sigma_p\in \gal\bigr(F(\zeta_r)/F\bigr)$ acting as $\sigma_p(\zeta_r)=\zeta_r^p$, and $L$ is the subfield of $F(\zeta_r)$ fixed by it, so that $F(\zeta_r)/L$ is cyclic of degree $s$. Finally we obtain an isomorphism 
$$\Q[\mathsf{f}_0,\tau] \xrightarrow{\sim} \bigl( F(\zeta_r)/L,\sigma_p,\pi \bigr)$$
given by $\tau \mapsto \zeta_r$ and $\mathsf{f}_0 \mapsto u_p$, with $u_p$ acting as $\sigma_p$ on $F(\zeta_r)$ and such that $u_p^s=\pi$. 
\end{proof}

\section{Representations arising from abelian varieties over $\F_p$}
\label{sec:repnsabvarFp}

\subsection{Semisimple representations defined over $\Q$}
\label{sec:ssQ}
Given a semisimple representation defined over $\Q$ it is not necessarly true that a subrepresentation is also defined over $\Q$. However lemma~\ref{lem1} shows that it is the case for subobjects of a certain type. Lemma~\ref{arithmlem} deals with the arithmetic of abelian semisimple representations defined over $\Q$.

\medskip
Let $L/K$ be a characteristic zero field extension and $V$ a finite dimensional $L$-vector space. A subgroup $H$ of $\aut_L(V)$ is {\em defined over $K$} if $\pchar(u)\in K[X]$ for all $u\in H$. If every element in $H$ is semisimple then this is equivalent to $\tr(u)\in K$ for all $u\in H$.

\begin{lem}
\label{lem1}
Let $f$ and $g$ be two semisimple elements in $\Endom_L(V)$ such that $\tr(g^n)\in K$ and $\tr(fg^n)\in K$ for all integers $n$. Let $W$ be a subspace of the form $W=\Ker P(g)$ for some $P\in K[X]$ and assume $f(W)\subseteq W$. Then $\tr(f_{\mid W})\in K$.
\end{lem}

\begin{proof}
Put $h=P(g)\in \Endom_L(V)$: it is a semisimple element such that $\tr(fh^n)\in K$ for all $n$, and $W=\Ker h$. Moreover, $\pchar(h)\in K[X]$ because $P \in K[X]$ and $\pchar (g) \in K[X]$ by assumption. Write $\pchar(h)(X)= X^m R(X)$ with $R(X)\in K[X]$ and $R(0)\neq 0$, which gives the decomposition $V=W\oplus \Ker R(h)$. Now $fR(h)$ stabilises $W$ with 
$$\bigl(fR(h)\bigr)_{\mid W}= R(0)f_{\mid W}$$ 
and vanishes on $\Ker R(h)$, so that
$$\tr\bigl(fR(h)\bigr) = R(0)\tr(f_{\mid W}).$$
Since $\tr\bigl(fR(h)\bigr) \in K$ and $R(0) \in K^{\times}$ we get $\tr(f_{\mid W})\in K$.
\end{proof}

\begin{cor}
\label{sbrpnscor}
Let $(\Delta,\rho)$ be a semisimple Weil representation defined over $\Q$. Let $u\in \im\rho$, $P\in \Q[X]$ such that $P(u)=0$, and $P=P_1P_2$ with $P_i$ coprime in $\Q[X]$. If $\Delta_i=\Ker P_i(u)$ is stable by the action of $\rho$ for $i=1,2$, then $\Delta=\Delta_1\oplus \Delta_2$ with each $\Delta_i$ defined over $\Q$.
\end{cor}

\begin{proof}
Apply lemma~\ref{lem1} to $K=\Q$, $f\in \im\rho$, $g=u$ and $W=\Delta_i$.
\end{proof}

\medskip
For $u\in \Endom_L(V)$ let $\spec u$ be the set of $u$'s eigenvalues in an algebraic closure $\overline{L}$ of $L$.

\begin{lem}
\label{arithmlem}
Let $u$ and $v$ be two semisimple commuting elements in $\aut_L(V)$ such that $\langle u,v\rangle$ is defined over $K$ and $\pchar(u)=\pmin(u)$. Then $\spec v \subset K(\spec u)$.
\end{lem}

\begin{proof}
The assumptions imply that $u$ and $v$ are codiagonalisable in $\overline{L}\otimes_L V$ and that $\pchar(u)$ has distinct roots. Let $d=\dim_L V$. By semisimplicity $\langle u,v\rangle$ is defined over $K$ if and only if $\tr(u^iv^j)\in K$ for all $0\leq i,j \leq d-1$. For $x_1,\dots,x_d \in \overline{L}$ let ${\rm V}(x_1,\dots,x_d)$ be the Vandermonde matrix
$${\rm V}(x_1,\dots,x_d) = 
\begin{pmatrix}
1         &  1         &  \dots    & 1         \\
x_1       &  x_2       &   \dots   & x_d       \\
x_1^2     &  x_2^2     &  \dots    & x_d^2     \\
\vdots    & \vdots     &           & \vdots     \\
x_1^{d-1} & x_2^{d-1}  &  \dots    & x_d^{d-1}  
\end{pmatrix}$$
Pick a basis of $\overline{L}\otimes_L V$ over which the matrices of $u$ and $v$ are respectively $\diag(\alpha_1, \dots,\alpha_d)$ and $\diag(\beta_1, \dots,\beta_d)$. Let $T$ be the matrix $\bigl(\tr(u^iv^j)\bigr)_{i,j} \in M_d(K)$. We have 
$${\rm V}(\alpha_1, \dots,\alpha_d)  {\rm V}(\beta_1, \dots,\beta_d)^t =T$$
where $M^t$ is the transpose of $M$. But ${\rm V}(\alpha_1, \dots,\alpha_d)$ is invertible by assumption, so 
$${\rm V}(\beta_1, \dots,\beta_d)^t = {\rm V}(\alpha_1, \dots,\alpha_d)^{-1} T$$
with $T\in M_d(K)$, which gives $\beta_r \in K(\alpha_1, \dots,\alpha_d)$ for all $1\leq r\leq d$.
\end{proof}

\subsection{Tame representations arising from Galois pairs}
\label{sec:galpairsrepns}
A semisimple tame Weil representation defined over $\Q$ is the direct sum of some distinguished subobjects, enjoying the same properties, that we call $\Q$-elementary (lemma~\ref{Q-elemdecomplem}). Each $\Q$-elementary object carries a specific arithmetic structure which enables us to attach some invariant to it. Requiring its dimension to be a multiple of this invariant leads to the notion of Tate dimension. A semisimple tame Weil representation defined over $\Q$ has Tate type if each of its $\Q$-elementary component has Tate dimension. This notion is part of the necessary and sufficient conditions stated in theorem~\ref{HTthm} for a tame $\bigl(\varphi,\gal(K/\Qp)\bigr)$-module to arise from a Galois pair. There is an obvious equivalent description for Weil representations.

\medskip
We fix some notation. Let $\Delta=(\Delta,\rho)$ be a semisimple tame Weil representation defined over $\Q$. Let $\phi_0$ be the image by $\rho$ of a lifting of the geometric Frobenius and $\theta$ a generator of $\rho(I)$. Then $\im \rho =\langle \phi_0,\theta \rangle$ with the relation $\phi_0\theta =\theta^p\phi_0$. The algebra $\Q[\phi_0,\theta]\subset \Endom_{K_0}(\Delta)$ spanned over $\Q$ by $\im\rho$ is finite dimensional over $\Q$ (because the representation is defined over $\Q$) and semisimple (because the representation is). As usual $\zeta_r\in \Qa$ is a primitive $r$th root of unity and $\Phi_r$ is the $r$th cyclotomic polynomial. The following lemma is a representation analogue of proposition~\ref{decompprop}.

\begin{lem}
\label{Q-elemdecomplem}
A semisimple tame Weil representation $\Delta$ defined over $\Q$ decomposes as 
$$\Delta = \bigoplus \Delta(r;\pi)$$
where $r$ is an integer prime to $p$, $\pi$ is a root of $\pmin(\phi_0^s)$ with $s=\ord\bigl(p\in(\Z/r\Z)^{\times}\bigr)$, and $\Delta(r;\pi)$ is a representation defined over $\Q$ enjoying the following properties:
\begin{itemize}
\item[(a)] $\pmin(\theta)=\Phi_r$ and $Z(\Q[\phi_0,\theta])\simeq L=$ a number field,
\item[(b)] $[F(\zeta_r):\Q]$ divides $\dim \Delta(r;\pi)$, with $F=\Q(\pi)$, 
\item[(c)] $\gal(F(\zeta_r)/F)$ contains an element $\sigma_p$ such that $\sigma_p(\zeta_r)= \zeta_r^p$, and
\item[(d)] $\Q[\phi_0,\theta] \simeq \bigl( F(\zeta_r)/L,\sigma_p,\pi \bigr)$ by $\theta \mapsto \zeta_r$, $\phi_0 \mapsto u_p$, with $u_p$ acting as $\sigma_p$ on $F(\zeta_r)$.
\end{itemize}
\end{lem}

\begin{proof}
The representation being semisimple and defined over $\Q$ we have $\pmin(\theta)=\prod_{r}\Phi_r$ with $r\mid \ord(\theta)$ for all $r$ and $\lcm(\text{all }r)=\ord(\theta)$. Since $r$ is prime to $p$ we see that $\phi_0$ stabilises each $\Ker\Phi_r(\theta)$, and by corollary~\ref{sbrpnscor} the representation restricted to $\Ker\Phi_r(\theta)$ is defined over $\Q$. So we may assume that $\pmin(\theta)=\Phi_r$ for some $r$. Let $s$ be the order of $p$ in $(\Z/r\Z)^{\times}$. Then we take $K_0=\Frac W(\F_{p^s})$ and $\phi =\phi_0^s$ is the smallest power of $\phi_0$ commuting with $\theta$. Again by corollary~\ref{sbrpnscor} we may assume that $\pmin(u)$ is irreducible in $\Q[X]$ for all $u$ in $Z(\Q[\phi_0,\theta])$, i.e. that the centre of $\Q[\phi_0,\theta]$ is a field. The element $\phi_0$ acts by conjugation $c_0 : u\mapsto \phi_0u\phi_0^{-1}$ on $\Q[\phi,\theta]$ and the same argument as in the proof of proposition~\ref{decompprop} shows that $\Q[\phi,\theta]$ is a field. 

Codiagonalise $\phi$ and $\theta$ in $\Qpa \otimes_{K_0}\Delta$ and let $\pi\in \Qa$ be such that $\pi$ and $\zeta_r$ are the eigenvalues, for $\phi$ and $\theta$ respectively, of some common eigenvector. Put $F=\Q(\pi)$. If $u\in \Q[\phi,\theta]$ writes as $u=P(\phi,\theta)$ for some $P\in \Q[X,Y]$ then $\xi =P(\pi,\zeta_r)\in \spec u$. Since $\pmin(u)$ is irreducible in $\Q[X]$ we have 
$$\pmin(u)=\pmin(\xi)$$
the right-hand side being the minimal polynomial over $\Q$ of an algebraic number. Pick a primitive element $\xi_0$ for $F(\zeta_r)/\Q$, i.e. such that $F(\zeta_r)=\Q(\xi_0)$. Let $P_0\in \Q[X,Y]$ be such that $P_0(\pi,\zeta_r)=\xi_0$ and put $u_0=P_0(\phi,\theta)$. Then $\pmin(u_0)=\pmin(\xi_0)$, so $\dim_{\Q}\Q[u_0]=[F(\zeta_r):\Q]$ which divides $\dim_{K_0}\Delta$, and $\Q[\phi,\theta]=\Q[u_0]$. The choice of $\xi_0$ yields a $\Q$-algebras isomorphism 
$$\Q[\phi,\theta] \xrightarrow{\sim} F(\zeta_r)$$ 
mapping $u_0$ to $\xi_0$ (thus $\phi \mapsto \pi$ and $\theta \mapsto \zeta_r$). By transport of structure the automorphism $c_0$ is carried to an element $\sigma_p\in \gal(F(\zeta_r)/F)$ acting as $\zeta_r \mapsto \zeta_r^p$. Hence the above isomorphism identifies $Z(\Q[\phi_0,\theta])$ with the subfield $L$ of $F(\zeta_r)$ fixed by $\sigma_p$. Finally we get an isomorphism 
$$\Q[\phi_0,\theta] \xrightarrow{\sim} \bigl( F(\zeta_r)/L,\sigma_p,\pi \bigr)$$
by sending $\phi_0$ to an element $u_p$ acting as $\sigma_p$ on $F(\zeta_r)$ and such that $u_p^s=\pi$.
\end{proof}

\medskip
Note that the representation $\Delta(r;\pi)$ is determined up to isomorphism by its dimension together with property (d) in lemma~\ref{Q-elemdecomplem}. Call such a $\Delta(r;\pi)$ {\em $\Q$-elementary} and $\Delta = \bigoplus \Delta(r;\pi)$ the $\Q$-elementary decomposition of $\Delta$.

\medskip
Now assume that the roots of $\pchar(\phi_0)$ are $p$-Weil numbers, hence $\pi$ is a $p^s$-Weil number. Put $n(1;\pi)=n(2;\pi)=1$, and for $r\geq 3$ recall that we have defined in section~\ref{sec:HTalg} the integer
$$n(r;\pi) \underset{\text{def}}{=} 
\begin{cases}
2\; \text{ if $L_v=\R$ and $\pi= p^{s/2}$ with $s$ even,} & \\
\lcm \biggl[ \ord \Bigl( \pi \in L_v^{\times} {\big /} N_{F(\zeta_r)_v/L_v}\bigl(F(\zeta_r)_v^{\times}\bigr) \Bigr)\text{, $v\mid \l$ in $L$ and $\l \mid r$}\biggr] \text{ otherwise.} & 
\end{cases}$$

\begin{defn}
\label{Tatetypedefn}
A $\Q$-elementary object $\Delta(r;\pi)$ as in lemma~\ref{Q-elemdecomplem} has {\em Tate dimension} if $[F(\zeta_r):\Q]n(r;\pi)$ divides $\dim_{K_0}\Delta(r;\pi)$. A semisimple tame Weil representation defined over $\Q$ has {\em Tate type} if each object of its $\Q$-elementary decomposition has Tate dimension. 
\end{defn}

\begin{ex}
\label{cyclicpairrepnex}
Let $(A_0,\langle \tau \rangle)$ be a cyclic tame Galois pair with $\Q[\mathsf{f}_0,\tau]\simeq \bigl( F(\zeta_r)/L,\sigma_p,\pi \bigr)$ as in proposition~\ref{decompprop}. Then $\Delta(A)$ is isomorphic to a $\Q$-elementary $\Delta(r;\pi)$ and proposition~\ref{embedprop} shows that it has Tate dimension. 
\end{ex}

\begin{rem}
\label{ladicrem}
Say that a Weil representation $\Delta$ is {\em $\l$-adically realisable} if there exists a compatible system of $\l$-adic Weil representations $(\Delta_{\l})_{\l}$, $\l$ running over all rational primes, such that $\Delta \simeq \Delta_p$ (in particular $\Delta$ is defined over $\Q$). Now let $\Delta$ be semisimple, defined over $\Q$, and tame. One checks that $\Delta$ has Tate type if and only if it is $\l$-adically realisable and $2\boldsymbol{\varphi}(r) \mid \dim\Delta(r;p^{s/2})$ for each $\Q$-elementary subobject $\Delta(r;p^{s/2})$ with $r\geq 3$ and $L$ totally real.
\end{rem}

\medskip
Let $K/\Qp$ be a finite tame Galois extension and $K_0$ the maximal unramified subfield of $K$. We have $\gal(K/\Qp)= I(K/\Qp)\rtimes \gal(K_0/\Qp)$. Let $D$ be a $\bigl(\varphi,\gal(K/\Qp)\bigr)$-module. Recall that the restriction of $\varphi$ on the sub-$\Qp$-vector space of $D$ consisting of elements fixed by $\gal(K_0/\Qp)$ is a $\Qp$-linear isomorphism $\varphi_0$. Consider the following conditions on $D$:
\begin{itemize}
\item[(1)] $\varphi_0$ acts semisimply and $\pchar(\varphi_0)$ is a $p$-Weil polynomial
\item[(2)] $\W(D)$ is defined over $\Q$ and has Tate type
\item[(3)] There exists a nondegenerate skew form $D\times D\rightarrow K_0\{-1\}$.
\end{itemize}
The skew form in $(3)$ is to be understood in the category of $\bigl(\varphi,\gal(K/\Qp)\bigr)$-modules: $\varphi_0$ is a $p$-similitude and $I(K/\Qp)$ acts by isometries. Note that the semisimplicity condition in $(1)$ implies that $\W(D)$ is semisimple (it is actually equivalent).

\begin{thm}
\label{HTthm}
Let $K/\Qp$ be a finite tame Galois extension and $D$ a $\bigl(\varphi,\gal(K/\Qp)\bigr)$-module. The following are equivalent:
\begin{itemize}
\item[(i)] There exists a Galois pair $(A_0,\langle\tau\rangle)$ for $K/\Qp$ such that $D\simeq D(A)$
\item[(ii)] $D$ satisfies conditions $(1)$, $(2)$, and $(3)$.
\end{itemize}
\end{thm}

In other words, a tame $\bigl(\varphi,\gal(K/\Qp)\bigr)$-module satisfying $(1)$ is ``geometric'' if and only if it also satisfies $(2)$ and $(3)$.

\begin{proof}
Let $D(A)$ be the $\bigl(\varphi,\gal(K/\Qp)\bigr)$-module associated to some tame Galois pair for $K/\Qp$. That $D(A)$ satisfies $(1)$ and is defined over $\Q$ is well-known. That it has Tate type follows from propositions~\ref{decompprop} and~\ref{embedprop}. That it satisfies $(3)$ follows from proposition~\ref{polgalprop}. 

\medskip 
Now let $D$ be a $\bigl(\varphi,\gal(K/\Qp)\bigr)$-module satisfying conditions $(1)$, $(2)$, and $(3)$. Clearly we may work with Weil representations, i.e. it is enough to show the existence of a tame Galois pair $(A_0,\langle\tau\rangle)$ for $K/\Qp$ such that $\W(D)=\Delta\simeq \Delta(A)$. Recall that $\phi_0=\W(\varphi_0)$. Let $\Delta = \bigoplus \Delta(r;\pi)$ be the $\Q$-elementary decomposition as in lemma~\ref{Q-elemdecomplem}. The $W$-equivariant nondegenerate skew form $\Delta \times \Delta \rightarrow K_0\{-1\}$ furnished by condition $(3)$ restricts to such a form on each $\Delta(r;\pi)$. In particular the restriction of $\phi_0$ is a symplectic $p$-similitude, which implies that $\pchar(\phi_0)$ remains a $p$-Weil polynomial on $\Delta(r;\pi)$. Indeed, all we have to check is that its valuation at $X^2-p$ stays even; but if it wasn't the case the determinant would be $-p^n$ for some $n$, a contradiction. Hence we may assume that $\Delta =\Delta(r;\pi)$ and $K_0=\Frac W(\F_{p^s})$ with $s=\ord\bigl(p\in(\Z/r\Z)^{\times}\bigr)$. 

Condition $(1)$ implies the existence of an abelian variety $A_0$ over $\F_p$ such that $\pchar(\mathsf{f}_0)=\pchar(\phi_0)$. If $r=1$ or $2$ we obviously take $\tau=1$ and $\tau=-1$ respectively, so we may assume $r\geq 3$. We have $\dim \Delta = 2\dim A$ and 
$$\pchar(\mathsf{f})=\pchar(\phi_0^s)=\pmin(\pi)^m$$ 
for some integer $m$. Put $F=\Q(\pi)$ and let $\delta$ be the invariant attached to the $p^s$-Weil number $\pi$ (see~\ref{sec:HT}). By lemmas~\ref{splitlem} and~\ref{Q-elemdecomplem} (b) we have the chain divisibility 
$$\delta \  \mid \  [F(\zeta_r):F] \  \mid \  m = \frac{\dim \Delta}{[F:\Q]}.$$ 
Since $\pmin(\pi)$ is irreducible $A$ is isogenous to the $(m/\delta)$th power of some simple abelian variety over $\F_{p^s}$ with Honda-Tate division algebra $D$ of degree $\delta$. Hence 
$$\Endom_{\F_{p^s}}^{\circ}(A) = M_{m/\delta}(D).$$
Condition $(2)$ means that $\Delta$ has Tate dimension, i.e. $[F(\zeta_r):F]n(r;\pi)$ divides $m$. Then proposition~\ref{embedprop} implies the existence of an embedding   
$$\bigl(F(\zeta_r)/L,\sigma_p,\pi\bigr) \hookrightarrow \Endom_{\F_{p^s}}^{\circ}(A)$$
sending $u_p$ to $\mathsf{f}_0$, with $u_p$ an element acting as $\sigma_p$ on $F(\zeta_r)$ and such that $u_p^s=\pi$. Letting $\tau$ be the image of $\zeta_r$ under this embedding we get
$$\bigl(F(\zeta_r)/L,\sigma_p,\pi\bigr) \simeq \Q[\mathsf{f}_0,\tau] \subseteq \Endom_{\F_{p^s}}^{\circ}(A).$$
Now take a maximal order in $\Endom_{\F_{p^s}}^{\circ}(A)$ containing both $\mathsf{f}_0$ and $\tau$. By Waterhouse's result~\cite{Wa} Thm.3.13, this maximal order occurs as an endomorphism ring (over $\F_{p^s}$) of some abelian variety (defined over $\F_p$ and isogenous to $A_0$). Replacing $A_0$ by it we have $\tau \in \aut_{\F_{p^s}}(A)$. Then $(A_0,\langle\tau\rangle)$ is a Galois pair for $K/\Qp$ with associated Weil representation $\Delta(A)$ isomorphic to $\Delta$.
\end{proof}

\subsection{Some examples in low dimension}
\label{sec:ex}
We keep the notations of section~\ref{sec:galpairsrepns}. 

\subsubsection{Elliptic curves}
\label{sec:ellip}
Let $\Delta$ be a $2$-dimensional Weil representation satisfying the conditions of theorem~\ref{HTthm} such that $\pchar(\theta)=\Phi_e$ with $\boldsymbol{\varphi}(e)=2$, i.e. $e\in \{ 3,4,6\}$. The order of $p$ in $(\Z/e\Z)^{\times}$ is either $1$ or $2$. Let $E/\F_p$ be an elliptic curve whose Frobenius has characteristic polynomial $\pchar(\phi_0)$. Recall that $E$ is ordinary if $\tr(\phi_0)$ is prime to $p$ and supersingular otherwise. 

\medskip
If $p\equiv 1\bmod e$ the representation is abelian and lemma~\ref{arithmlem} implies $\Endom_{\F_{p}}^{\circ}(E)=\Q(\zeta_e)$. This is a quadratic extension of $\Q$ unramified at $p$, so $\tr(\phi_0)$ is prime to $p$ and $E$ is ordinary. Then $\tau =\zeta_e$ or $\zeta_e^{-1}$ (yielding nonisomorphic representations).  

If $p\equiv -1\bmod e$ the relation $\phi_0\theta =\theta^{-1}\phi_0$ implies $\tr(\phi_0)=0$. Thus $E$ is supersingular with $\zeta_e\in \Endom_{\F_{p^2}}^{\circ}(E)= D_{p,\infty}$. Then we may take $\tau =\zeta_e$ (taking $\zeta_e^{-1}$ yields an isomorphic representation).  

\medskip
A complete description of the Weil representations arising from elliptic curves over $\F_p$ when $p>3$ can be found in~\cite{Vo} sections 2.1 and 3.1.

\subsubsection{The product of two supersingular elliptic curves}
\label{sec:prodellip}
Let $\Delta$ be a $4$-dimensional Weil representation satisfying the conditions of theorem~\ref{HTthm} such that $\pchar(\theta)=\Phi_8$ (so $p$ is odd) and $\pchar(\phi_0)(X)=(X^2+p)^2$. Let $E/\F_p$ be a supersingular elliptic curve whose Frobenius has characteristic polynomial $X^2+p$. We have $\Endom_{\F_{p}}^{\circ}(E)=\Q(\sqrt{-p})$ and $\Endom_{\F_{p^2}}^{\circ}(E)= D_{p,\infty}$. Take $A_0=E\times E$ with $\mathsf{f}_0 = \diag ({\rm Frob}_E,{\rm Frob}_E) \in M_2\bigl(\Q(\sqrt{-p})\bigr)$.

\medskip
If $p\equiv 1\bmod 8$ then by lemma~\ref{arithmlem} we would have $\Q(\sqrt{-p})\subset \Q(\zeta_8)$, which is impossible since $\Q(\zeta_8)$ is unramified at $p$. 

If $p\equiv 3\bmod 8$ then $\zeta_4\in D_{p,\infty}$ and we may take
$\tau = 
\begin{pmatrix}
0 & -\frac{1}{2}+\frac{1}{2}\zeta_4 \\
1-\zeta_4 & 0
\end{pmatrix}.$

If $p\equiv 5\bmod 8$ then $\xi\in D_{p,\infty}$ with $\xi^2+2=0$ and we may take
$\tau = 
\begin{pmatrix}
\xi^{-1} & \xi^{-1} \\
-\xi^{-1} & \xi^{-1}
\end{pmatrix}.$

If $p\equiv -1\bmod 8$ then $\zeta_4\in D_{p,\infty}$ and we may take
$\tau = 
\begin{pmatrix}
0 & \frac{1}{2}+\frac{1}{2}\zeta_4 \\
1+\zeta_4 & 0
\end{pmatrix}.$

\section{Representations arising from abelian varieties over $\Qp$}
\label{sec:repabvar}

\subsection{Polarisations and Rosati involutions}
\label{sec:pol}
Let $k$ be a finite extension of $\F_p$ and $A/k$ an abelian variety. Fix a polarisation $\lambda =\lambda_{{\mathcal L}}: A \rightarrow A^{\vee}$ induced by an ample invertible sheaf ${\mathcal L}$ on $A$ (\cite{Mi} Rmk.13.2). The associated Rosati involution $\dag$ on $\Endom_k^{\circ}(A)$ is given by 
$$\alpha \mapsto \alpha^{\dag} = \lambda^{-1}\alpha^{\vee}\lambda.$$ 
It is positive definite, i.e. the bilinear form $(\alpha,\beta)\mapsto\tr(\alpha^{\dag}\beta)$ is. If $\alpha^{\dag}=\alpha$ then all the roots of its reduced characteristic polynomial $\prd(\alpha)$ lie in $\R$. Write $\Endom_k^{\circ}(A)^{\dag}$ for the sub-$\Q$-vector space of $\Endom_k^{\circ}(A)$ consisting of elements fixed by $\dag$ and define
$$\pol(A) \underset{\text{def}}{=} \{ \alpha \in \Endom_k^{\circ}(A)^{\dag} \mid \roots \bigl( \prd(\alpha)\bigr) \subset \R_{>0} \}.$$
All polarisations in $\Hom_k^{\circ}(A,A^{\vee})$, i.e. elements $\mu$ such that $m\mu$ is a polarisation on $A$ for some positive integer $m$, have the form 
$$\mu = \lambda \,\alpha  \quad \text{ for some }  \alpha \in \pol(A).$$
Indeed, an isogeny $\mu : A \rightarrow A^{\vee}$ is induced by an invertible sheaf ${\mathcal L}$ if and only if it is symmetric, i.e. $\lambda^{-1}\mu$ is fixed by $\dag$ (\cite{Mi} Prop.17.2), and then ${\mathcal L}$ is ample if and only if all the roots of $\prd(\lambda^{-1}\mu)$ are positive (\cite{Mu} III.16, Thm.p.155).

\medskip
If $D=D(A)$ then $D(A^{\vee})=D^*\{-1\}$ and $\lambda$ yields an isomorphism of $\varphi$-modules
$$\delta =\D\bigl(\lambda(p)\bigr) : D^*\{-1\} \xrightarrow{\sim} D \quad \text{ such that } \quad \delta^*\{-1\}= - \delta$$
(under the canonical identification $\bigl(D^*\{-1\}\bigr)^*\{-1\}\simeq D$). Indeed, a symmetric morphism such as $\lambda$ (i.e. $\lambda^{\vee}=\lambda$ under ${A^{\vee}}^{\vee}\simeq A$) becomes antisymmetric when passing to $p$-divisible groups, since $A^{\vee}(p)$ is the Cartier dual of $A(p)$. Thus Tate's isomorphism restricts to 
$$\Qp \otimes_{\Q} \Hom_k^{\text{s}} (A,A^{\vee}) \simeq \Hom_{\varphi}^{\text{a}}(D^*\{-1\},D)$$
where the superscripts ``s'' and ``a'' stand respectively for ``symmetric'' in $\Hom_k^{\circ} (A,A^{\vee})$ and ``antisymmetric''. Of course, an antisymmetric isomorphism of $\varphi$-modules such as $\delta$ defines a nondegenerate alternating form in the category of $\varphi$-modules 
$$b: D \times D \rightarrow K_0\{-1\} \quad \text{ with }K_0=\Frac W(k)$$ 
(i.e. $b(\varphi x,\varphi y)=q b(x,y)$ if $k=\F_q$), and vice versa.

\medskip
The Rosati involution extends by linearity (or continuity) to $\Qp\otimes_{\Q}\Endom_k^{\circ}(A)$ and yields an involution also denoted $\dag$ on $\Endom_{\varphi}(D)$, given by $\alpha \mapsto \delta\alpha^*\{-1\}\delta^{-1}$. Again write $\Endom_{\varphi}(D)^{\dag}$ for the sub-$\Qp$-vector space of elements fixed by $\dag$. We have the following elementary result.

\begin{lem}
\label{poldenselem}
The image of $\pol(A)$ in $\Endom_{\varphi}(D)^{\dag}$ is dense.
\end{lem}

\begin{proof}
If $x\in \Endom_k^{\circ}(A)^{\dag}$ then $x+p^n\id \in \pol(A)$ for a sufficiently large integer $n$. Therefore $\pol(A)$ is dense in $\Endom_k^{\circ}(A)^{\dag}$ with respect to the $p$-adic topology. The lemma then follows from Tate's theorem $\Qp\otimes_{\Q}\Endom_k^{\circ}(A)^{\dag} \simeq \Endom_{\varphi}(D)^{\dag}$.  
\end{proof}

\medskip
Let $E$ be a finite dimensional semisimple algebra over $\Qp$ and $\dag$ any involution on $E$ fixing $\Qp$. The group $E^{\times}$ of invertible elements acts on the $\Qp$-vector space $E^{\dag}$ by $\gamma \mapsto \alpha \gamma \alpha^{\dag}$, $\gamma \in  E^{\dag}$, $\alpha \in E^{\times}$. The proof of the following well-known result is left to the reader. 

\begin{lem}
\label{openorblem}
The orbit of an invertible element in $E^{\dag}$ under the above action of $E^{\times}$ is open.
\end{lem}

\subsection{Polarisable Galois pairs}
\label{sec:polgalpairs}
Let $(A_0,\Gamma)$ be a Galois pair for $K/\Qp$ with 
$$\nu : \gal(K/\Qp) \hookrightarrow \ist_{\F_{p^s}}(A).$$ 
Put $\Gamma^{\vee} = \{ \gamma^{\vee}, \gamma \in \Gamma \} \subseteq \aut_{\F_{p^s}}(A^{\vee})$ and write $\iota : \ist_{\F_{p^s}}(A) \xrightarrow{\sim} \ist_{\F_{p^s}}(A^{\vee})$ for the group isomorphism given by $\iota(\xi)=(\xi^{\vee})^{-1}$. Then $(A_0^{\vee},\Gamma^{\vee})$ is a Galois pair for $K/\Qp$ with
$$\nu^{\prime} = \iota \circ \nu : \gal(K/\Qp) \hookrightarrow \ist_{\F_{p^s}}(A^{\vee}).$$
The $(\varphi,\gal(K/\Qp))$-module $D(A^{\vee})$ obtained from $\nu^{\prime}$ is the twisted dual $D(A)^*\{-1\}$ of the $(\varphi,\gal(K/\Qp))$-module $D(A)$ obtained from $\nu$, since the dual action of an element in $I(K/\Qp)$ is the inverse of the dual of its image in $\aut_{\varphi}\bigl(D(A)\bigr)$.

\medskip
Let $\lambda_0 : A_0 \rightarrow A_0^{\vee}$ be an $\F_p$-isogeny and $\lambda = \lambda_0\times_{\F_p}\F_{p^s}: A \rightarrow A^{\vee}$. Then $\delta =\D\bigl(\lambda(p)\bigr) : D(A^{\vee}) \xrightarrow{\sim} D(A)$ is an isomorphism of $\bigl(\varphi,\gal(K/\Qp)\bigr)$-modules if and only if $\delta \D\bigl(\gamma^{\vee}(p)^{-1}\bigr) = \D\bigl(\gamma(p)\bigr)\delta$ for all $\gamma \in \Gamma$, that is,
$$\gamma^{\vee} \lambda \gamma = \lambda \quad \text{ for all } \gamma \in \Gamma.$$
When $\lambda_0$ is a polarisation this means $\Gamma \subseteq \aut(A,\lambda)$ (compare with lemma~\ref{pgrpollem}).

\begin{defn}
\label{polgaldefn}
The Galois pair $(A_0,\Gamma)$ is {\em polarisable} if there exists an $\F_p$-polarisation $\lambda_0$ on $A_0$ such that $\Gamma \subseteq \aut(A,\lambda)$.
\end{defn}

\begin{prop}
\label{polgalprop}
Any tame Galois pair is polarisable.
\end{prop}

\begin{proof}
Let $(A_0,\langle \tau \rangle)$ be a tame Galois pair with $\tau \in \aut_{\F_{p^s}}(A)$. By proposition~\ref{decompprop} we may assume it to be cyclic, so $\Endom_{\F_{p^s}}^{\circ}(A)$ is central simple over $F$. Pick any $\F_p$-polarisation $\lambda_0$ on $A_0$ with associated Rosati involution $\dag$. Define $\mathcal{I}_{0}^{+}$ to be the set of positive definite involutions on $\Endom_{\F_{p^s}}^{\circ}(A)$ coinciding with $\dag$ on $F$ (thus having same kind), having same type as $\dag$ (in case it is of the first kind), and stabilising $\Endom_{\F_{p}}^{\circ}(A_0)$. We have a map 
$$\begin{cases}
\pol(A_0) & \rightarrow  \quad  \mathcal{I}_{0}^{+}  \\
\alpha    & \mapsto      \quad  I_{\alpha} : x \mapsto \alpha^{-1}x^{\dag}\alpha \text{ for all } x \in \Endom_{\F_{p^s}}^{\circ}(A)
\end{cases}$$
The involution $I_{\alpha}$ is the Rosati involution associated to the $\F_p$-polarisation $\lambda_0\alpha$. This map is surjective by the Skolem-Noether theorem. Now put 
$$\pol(A_0,\tau) \underset{\text{def}}{=} \{ \alpha \in \pol(A_0) \mid \tau \in \aut(A,\lambda\alpha) \}.$$
We want to show that $\pol(A_0,\tau)$ is not empty, since for $\alpha$ in this set $\lambda_0\alpha$ is an appropriate polarisation. Put
$$\mathcal{I}_{0,\tau}^{+} \underset{\text{def}}{=} \{ \text{ involutions }\ddag \in \mathcal{I}_{0}^{+} \mid \tau^{\ddag}\tau = 1 \}.$$
Obviously the above map carries $\pol(A_0,\tau)$ onto $\mathcal{I}_{0,\tau}^{+}$, hence it is enough to show that $\mathcal{I}_{0,\tau}^{+}$ is nonempty. Define a positive definite involution ${}^{\prime}$ on $\Q[\mathsf{f}_0,\tau] \subseteq \Endom^{\circ}_{\F_{p^s}}(A)$ by setting
$$\mathsf{f}_0^{\prime} = p\mathsf{f}_0^{-1} \quad \text{ and } \quad \tau^{\prime} = \tau^{-1}.$$
Clearly $\dag$ and ${}^{\prime}$ have the same restriction to $F$. We claim that by Theorem 4.14 of~\cite{BI} the involution ${}^{\prime}$ can be extended to an involution on $\Endom^{\circ}_{\F_{p^s}}(A)$ belonging to $\mathcal{I}_{0,\tau}^{+}$. Indeed, since $\mathsf{f}_0^{\prime} = p\mathsf{f}_0^{-1}$ any extension of ${}^{\prime}$ stabilises $\Endom_{\F_{p}}^{\circ}(A_0)$, and since it is positive definite it can be extended to a positive definite one on $\Endom^{\circ}_{\F_{p^s}}(A)$. It remains to check that, in case ${}^{\prime}$ is of the first kind (and thus $\dag$ as well), it can be extended to an involution of the same type. So assume ${}^{\prime}$ to be of the first kind. Then $s$ must be even and $\mathsf{f}=\pm p^{s/2}$, so $F=\Q$, $\Q[\mathsf{f}_0,\tau]\simeq \bigl(\Q(\zeta_r)/L,\sigma_p,\pm p^{s/2}\bigr)$ with the notations of proposition~\ref{decompprop}, and
$$\Endom^{\circ}_{\F_{p^s}}(A) = M_n ( D_{p,\infty} ) \quad \text { for some integer } n.$$
The involution $\dag$ on $M_n(D_{p,\infty})$ has symplectic type. On the other hand ${}^{\prime}$ is symplectic if $\mathsf{f} = - p^{s/2}$ and orthogonal if $\mathsf{f} = p^{s/2}$. According to Theorem 4.14 of~\cite{BI} we need to check in the latter case that the degree of the centraliser $C$ of $\bigl( \Q(\zeta_r)/L ,\sigma_p, p^{s/2} \bigr)$ in $M_n ( D_{p,\infty} )$ is even, for then ${}^{\prime}$ can be extended to a symplectic involution. We have
$$\deg_L C = 2n / \boldsymbol{\varphi}(r).$$
Now $L$ is totally real of degree $\boldsymbol{\varphi}(r)/s$ over $\Q$ and $p^{s/2}$ is positive, so  
$$\R\otimes_{\Q} \bigl( \Q(\zeta_r)/L ,\sigma_p, p^{s/2} \bigr) \simeq \Bigl(  \R\otimes_{L} \bigl( \Q(\zeta_r)/L ,\sigma_p, p^{s/2} \bigr) \Bigr)^{\boldsymbol{\varphi}(r)/s}  \simeq  M_s(\R)^{\boldsymbol{\varphi}(r)/s}.$$
Therefore there is an embedding 
$$M_s(\R)^{\boldsymbol{\varphi}(r)/s} \hookrightarrow M_n(\H) \simeq \R\otimes_{\Q}M_n ( D_{p,\infty} )$$
which implies that $\frac{\boldsymbol{\varphi}(r)}{s}s = \boldsymbol{\varphi}(r)$ divides $n$. Thus $\deg_L C$ is even.
\end{proof}

\subsection{Lifting polarisations}
\label{sec:geomsp}
Let $(A_0,\langle \tau \rangle)$ be a tame Galois pair for $K/\Qp$ with $\tau \in \aut_{\F_{p^s}}(A)$. Proposition~\ref{polgalprop} allows us to pick an $\F_p$-polarisation $\mu_0$ on $A_0$ with associated Rosati involution $\dag$ such that $\tau \in \aut(A,\mu)$, which also writes as 
$$\tau^{\dag}\tau=1.$$ 
These data will be fixed throughout the rest of this section, as well as the notations $\Endom_{\F_p,\tau}^{\circ}(A_0)= \{\, x_0 \in \Endom_{\F_p}^{\circ}(A_0) \mid x\tau=\tau x \, \}$, $D_0=D(A_0)$, $D=D(A)$, $K_0=\Frac W(\F_{p^s})$, superscripts ``s'' and ``a'' for ``symmetric'' and ``antisymmetric'', and so on.

\medskip
The relation $\tau^{\dag}\tau=1$ implies that $\Endom_{\F_p,\tau}^{\circ}(A_0)$ and $\Endom_{\varphi,\gal(K/\Qp)}(D)$ are stable by $\dag$. Put
$$\Hom_{\F_p,\tau}^{\text{s}} (A_0,A_0^{\vee}) \underset{\text{def}}{=} \{ \, f \in \Hom_{\F_p}^{\text{s}} (A_0,A_0^{\vee}) \mid f\tau^{-1}=\tau^{\vee}f \, \}.$$
The isomorphism $\Endom_{\F_p}^{\circ}(A_0)^{\dag} \xrightarrow{\sim} \Hom_{\F_p}^{\text{s}}(A_0,A_0^{\vee})$ given by $x\mapsto \mu_0x$ carries $\Endom_{\F_p,\tau}^{\circ}(A_0)^{\dag}$ into $\Hom_{\F_p,\tau}^{\text{s}}(A_0,A_0^{\vee})$, again because $\tau^{\dag}\tau=1$. By Tate's theorem as stated in remark~\ref{Tatewithdescentrmk} we have a commutative diagram
$$\xymatrix{
\Endom_{\F_p,\tau}^{\circ}(A_0)^{\dag} \ar@{^{(}->}[r]_{\text{can}} \ar[d]^{\wr}  &  \Qp \otimes_{\Q}\Endom_{\F_p,\tau}^{\circ}(A_0)^{\dag} \ar[r]^{\sim} \ar[d]^{\wr}  &  \Endom_{\varphi,\gal(K/\Qp)}(D)^{\dag}  \ar[d]^{\wr}   \\
\Hom_{\F_p,\tau}^{\text{s}}(A_0,A_0^{\vee}) \ar@{^{(}->}[r]_{\text{can}}  &  \Qp \otimes_{\Q}\Hom_{\F_p,\tau}^{\text{s}}(A_0,A_0^{\vee}) \ar[r]^{\sim}  &\Hom_{\varphi,\gal(K/\Qp)}^{\text{a}}(D^*\{-1\},D)  }$$  
where the right-hand side vertical isomorphism is given by $\alpha\mapsto \alpha\nu_0$ with $\nu_0 = \D\bigl(\mu_0(p)\bigr)$. The following lemma has been pointed out by O. B\"ultel.

\begin{lem}
\label{oliverlem}
Assume there is an antisymmetric isomorphism $\delta : D^*\{-1\} \xrightarrow{\sim} D$ of $\bigl(\varphi,\gal(K/\Qp)\bigr)$-modules. Then there exists an $\F_p$-polarisation $\lambda_0$ on $A_0$ and an element $\alpha$ in $\aut_{\varphi,\gal(K/\Qp)}(D)$ such that $\D\bigl(\lambda(p)\bigr) = \alpha \delta \alpha^*$ and $\tau \in \aut(A,\lambda)$. 
\end{lem}

\begin{proof}
Write $\delta=\gamma\nu_0$ with $\gamma \in \aut_{\varphi,\gal(K/\Qp)}(D)^{\dag}$. Recall that $\pol(A_0,\tau)$ is the set of $\alpha \in \pol(A_0)$ such that $\tau \in \aut(A,\mu_0\alpha)$, so that the map $\beta \mapsto \mu_0\beta$  identifies $\pol(A_0,\tau)$ with the set of polarisations $\lambda_0$ on $A_0$ (up to a positive integer) such that $\tau \in \aut(A,\lambda)$. Now since $\tau^{\dag}=\tau^{-1}$ we have 
$$\pol(A_0,\tau)= \pol(A_0) \cap \Endom_{\F_p,\tau}^{\circ}(A_0) \subseteq \Endom_{\F_p,\tau}^{\circ}(A_0)^{\dag}.$$ 
The same argument as in lemma~\ref{poldenselem} shows that the image of $\pol(A_0,\tau)$ is dense in $\Endom_{\varphi,\gal(K/\Qp)}(D)^{\dag}$. Also, the group $\aut_{\varphi,\gal(K/\Qp)}(D)$ acts on $\Endom_{\varphi,\gal(K/\Qp)}(D)^{\dag}$ by $\xi \mapsto \alpha\xi\alpha^{\dag}$, $\xi \in \Endom_{\varphi,\gal(K/\Qp)}(D)^{\dag}$, $\alpha \in \aut_{\varphi,\gal(K/\Qp)}(D)$, and by lemma~\ref{openorblem} the orbit of $\gamma$ is open. These two facts imply the existence of an $\alpha \in \aut_{\varphi,\gal(K/\Qp)}(D)$ such that 
$$\alpha \gamma \alpha^{\dag} = \D\bigl(\beta(p)\bigr) \quad \text{ for some } \beta \in \pol(A_0,\tau).$$
Then $\lambda_0 = \mu_0\beta$ is the required polarisation since $\D\bigl(\lambda(p)\bigr) =  \alpha\gamma\alpha^{\dag}\nu_0 = \alpha\gamma (\nu_0\alpha^*\nu_0^{-1})\nu_0 = \alpha \delta \alpha^*$. 
\end{proof}

\medskip
Let $(D,\Fil)$ be a filtered $\bigl(\varphi,\gal(K/\Qp)\bigr)$-module of Hodge-Tate type $(0,1)$. A nondegenerate skew form on $(D,\Fil)$ is one such form $D\times D\rightarrow K_0\{-1\}$ that is a morphism of $\bigl(\varphi,\gal(K/\Qp)\bigr)$-modules and for which $\Fil^1D_K$ is totally isotropic (after extending the scalars to $K$). Equivalently it is given by an antisymmetric isomorphism of filtered $\bigl(\varphi,\gal(K/\Qp)\bigr)$-modules 
$$\delta : D^*\{-1\} \xrightarrow{\sim} D$$ 
since $D^*\{-1\}$ has Hodge-Tate type $(0,1)$ with $\Fil^1(D^*\{-1\})_K = (\Fil^1D_K)^{\perp}$. Now for $D=D(A)$, say that a polarisation $\lambda_0$ on $A_0$ lifts to $(D,\Fil)$ if it induces a (necessarly antisymmetric) isomorphism of filtered $\bigl(\varphi,\gal(K/\Qp)\bigr)$-modules $\D\bigl(\lambda(p)\bigr) : D^*\{-1\} \xrightarrow{\sim} D$.

\begin{prop}
\label{polprop}
Let $(A_0,\langle \tau \rangle)$ be a tame Galois pair for $K/\Qp$ and $D=D(A)$ the associated $\bigl(\varphi,\gal(K/\Qp)\bigr)$-module. Let $\Fil = (\Fil^iD_K)_{i\in \Z}$ be a Hodge-Tate $(0,1)$ filtration on $D_K$ stable by $\gal(K/\Qp)$. Assume there is a nondegenerate skew form on $(D,\Fil)$. Then there exists
\begin{itemize}
\item[(a)] a filtration $\Fil^{\prime}$ on $D_K$ such that $(D,\Fil)\simeq (D,\Fil^{\prime})$, and
\item[(b)] a polarisation $\lambda_0$ on $A_0$ lifting to $(D,\Fil^{\prime})$ and such that $\tau \in \aut(A,\lambda)$. 
\end{itemize}
\end{prop}

In other words, under the above assumption one can always find an object in the isomorphism class of $(D,\Fil)$ on which some polarisation on $A_0$, for which $\tau$ is a polarised automorphism, lifts.

\begin{proof}
The given nondegenerate skew form induces an antisymmetric isomorphism of filtered $\bigl(\varphi,\gal(K/\Qp)\bigr)$-modules $\delta : D^*\{-1\} \xrightarrow{\sim} D$. By lemma~\ref{oliverlem} there is a polarisation $\lambda_0$ on $A_0$ and an $\alpha \in \aut_{\varphi,\gal(K/\Qp)}(D)$ such that $\D\bigl(\lambda(p)\bigr) = \alpha \delta \alpha^*$ and $\tau \in \aut(A,\lambda)$. Put 
$$\Fil^{\prime} = (\alpha_K\Fil^iD_K)_{i\in \Z}.$$
Obviously $\alpha$ induces an isomorphism $(D,\Fil) \xrightarrow{\sim} (D,\Fil^{\prime})$ of filtered $\bigl(\varphi,\gal(K/\Qp)\bigr)$-modules. Furthermore, we have  
$$(\alpha\delta\alpha^*)_K(\alpha_K\Fil^1D_K)^{\perp} = \alpha_K \delta_K \bigl( \alpha^*_K (\alpha_K\Fil^1D_K)^{\perp} \bigr) = \alpha_K \delta_K (\Fil^1D_K)^{\perp} \subseteq \alpha_K \Fil^1D_K$$
so $\alpha\delta\alpha^* : D^*\{-1\} \xrightarrow{\sim} D$ is an isomorphism of $\bigl(\varphi,\gal(K/\Qp)\bigr)$-modules preserving the filtration $\Fil^{\prime}$. 
\end{proof}

\begin{rem}
We actually do not need the Hodge-Tate type to be $(0,1)$. However, the existence of an antisymmetric isomorphism of filtered modules $\delta : D^*\{-1\} \xrightarrow{\sim} D$ forces the Hodge-Tate type to be $(-n,n+1)$ for some nonnegative integer $n$.
\end{rem}

\subsection{The main theorem}
\label{sec:thm}
Let $K/\Qp$ be a finite tame Galois extension and $(D,\Fil)$ a filtered $\bigl(\varphi,\gal(K/\Qp)\bigr)$-module. As in section~\ref{sec:galpairsrepns} write $\varphi_0$ for the $\Qp$-linear restriction of $\varphi$ on $D^{\gal(K_0/\Qp)}$ with $K_0$ be the maximal unramified subfield of $K$. Consider the following conditions on $(D,\Fil)$:
\begin{itemize}
\item[(1)] $\varphi_0$ acts semisimply and $\pchar(\varphi_0)$ is a $p$-Weil polynomial
\item[(2)] $\W(D)$ is defined over $\Q$ and has Tate type
\item[(3)] There exists a nondegenerate skew form $(D,\Fil)\times (D,\Fil)\rightarrow K_0\{-1\}$
\item[(4)] It has Hodge-Tate type $(0,1)$.
\end{itemize}
When $(D,\Fil)\simeq\DcrK(V)$ condition $(3)$ means that there exists a $G$-equivariant symplectic form $V\times V\rightarrow \Qp(1)$. Hence ${\bigwedge}^{2d}V=\Qp(d)$ where $\dim V = 2d$, that is, the determinant on $V$ is the the $d$-th power of the $p$-adic cyclotomic character.

\begin{thm}
\label{mainthm}
Let $p\neq 2$. Let $V$ be a $p$-adic representation of $G$ that becomes crystalline over a finite tame Galois extension $K/\Qp$. The following are equivalent:
\begin{itemize}
\item[(i)] There exists an abelian variety ${\mathcal A}_0$ over $\Qp$ such that $V\simeq V_p({\mathcal A}_0)$
\item[(ii)] {\em $\DcrK(V)$} satisfies conditions $(1)$, $(2)$, $(3)$, and $(4)$.
\end{itemize}
\end{thm}

\begin{rem} When $V$ is $2$-dimensional conditions $(2)$ and $(3)$ may be replaced respectively by the weaker 
\begin{itemize}
\item[$(2^{\prime})$] $\W(D)$ is defined over $\Q$
\item[$(3^{\prime})$] ${\bigwedge}^2 (D,\Fil) =  K_0\{-1\}$  
\end{itemize}
the latter meaning ${\bigwedge}^{2}V=\Qp(1)$, see \cite{Vo}~Thm.5.1. An explicit list of the filtered $(\varphi,G)$-modules arising from elliptic curves over $\Qp$ when $p>3$ is given in~\cite{Vo} section~2.2.
\end{rem}

\begin{proof}
Let the abelian variety ${\mathcal A}_0/\Qp$ have good reduction over the tame Galois extension $K/\Qp$ and $(D,\Fil)=\DcrK\bigl(V_p({\mathcal A}_0)\bigr)$. By theorem~\ref{HTthm} the $\bigl(\varphi,\gal(K/\Qp)\bigr)$-module $D$ satisfies conditions $(1)$ and $(2)$. That $(D,\Fil)$ satisfies $(3)$ follows from the existence of a polarisation on ${\mathcal A}_0$ and lemma~\ref{pgrpollem}. That the filtration satisfies $(4)$ is well-known. 

\medskip
Now let $(D,\Fil)=\DcrK(V)$ be a (weakly admissible) filtered $\bigl(\varphi,\gal(K/\Qp)\bigr)$-module satisfying $(1)$, $(2)$, $(3)$, and $(4)$. We may assume that $K$ is the Galois closure of a totally ramified tame extension $L$ of minimal degree over which $V$ becomes crystalline. Let $\F_{p^s}$ be the residue field of $K$ and $K_0=\Frac W(\F_{p^s})$. By theorem~\ref{HTthm} conditions $(1)$, $(2)$, and $(3)$ imply the existence of a Galois pair $(A_0,\langle\tau\rangle)$ for $K/\Qp$ with associated $\bigl(\varphi,\gal(K/\Qp)\bigr)$-module $D(A)$ isomorphic to $D$. Hence we may assume that $D=D(A)$. By proposition~\ref{polprop} condition $(3)$ implies the existence of an $\F_p$-polarisation $\lambda_0$ on $A_0$ lifting to $(D,\Fil)$, possibly after replacing the filtration by an isomorphic one. By Breuil's theorem (\cite{Br} Thm.5.3.2) condition $(4)$ and $p\neq 2$ imply the existence of a $p$-divisible group ${\mathcal G}$ over the ring of integers $O_L$ of $L$ such that $V \simeq V_p({\mathcal G})$ as $G_L$-modules. Also, by a result of Raynaud~\cite{Ra} every $G_L$-stable lattice in $V$ comes from a $p$-divisible group as well. Replacing the Galois pair $(A_0,\langle\tau\rangle)$ by an isogenous one if necessary, we may assume that $A_0(p)$ lifts to such a $p$-divisible group over $O_L$. Let $T$ be a $G$-stable lattice in $V$ such that
$$T=T_p({\mathcal G}) \quad \text{with} \quad {\mathcal G}\times_{O_L}\F_p \simeq A_0(p).$$
The polarisation $\lambda_0$ induces an antisymmetric isomorphism of $G$-modules $\xi : V \xrightarrow{\sim} V^*(1)$. Multiplying $\lambda_0$ by a suitable power of $p$ we may assume that the restriction of $\xi$ on $T$ yields an injection $\xi : T \hookrightarrow T^*(1)$. By Tate's full faithfulness theorem~\cite{Ta2} we have 
$$\Hom_{\text{$p$-div/$O_L$}}\bigl({\mathcal G},{\mathcal G}^D\bigr) \simeq \Hom_{\Z_p[G_L]}\bigl(T,T^*(1)\bigr)$$
where ${\mathcal G}^D$ is the Cartier dual of ${\mathcal G}$, so that $T_p\bigl({\mathcal G}^D\bigr)=T^*(1)$. Therefore the isogeny $\lambda_0(p): A_0(p) \rightarrow A_0^{\vee}(p)= A_0(p)^D$ induced by $\lambda_0$ lifts to a quasipolarisation 
$$\Lambda(p) : {\mathcal G} \rightarrow {\mathcal G}^D.$$
By the Serre-Tate theory of liftings together with Grothendieck's theorem on algebraisation of formal schemes (\cite{Gr1} 5.4.5), the data $\bigl(A_0,\lambda_0,{\mathcal G},\Lambda(p)\bigr)$ defines an abelian scheme over $O_L$ lifting $A_0$. Its generic fibre base-changed to $K$ is an abelian variety ${\mathcal A}/K$ with special fibre $A/\F_{p^s}$ and $T\simeq T_p({\mathcal A})$ as $G_K$-modules. The action of $G_K$ extends to an action of $G$ on $T_p({\mathcal A})$ such that $G_L$ acts naturally (${\mathcal A}$ is defined over $L$) and $G_{K_0}$ induces an action of $I(K/\Qp)$ coming from $\langle\tau\rangle \subseteq \aut_{\F_{p^s}}(A)$. Therefore, by theorem~\ref{descthm} together with the subsequent comments, ${\mathcal A}$ is defined over $\Qp$, say ${\mathcal A} \simeq {\mathcal A}_0 \times_{\Qp}K$, and $V\simeq V_p({\mathcal A}_0)$ as $G$-modules.
\end{proof}

\begin{cor}
\label{cryscor}
Let $p\neq 2$. A crystalline $p$-adic representation of $G$ is isomorphic to the Tate module of an abelian variety over $\Qp$ if and only if its associated filtered module satisfies conditions $(1)$, $(3)$, and $(4)$ of theorem~\ref{mainthm}.
\end{cor}

\begin{cor}
\label{maincor}
Let $d$ be a positive integer and assume $p>2d +1$. A $2d$-dimensional potentially crystalline $p$-adic representation of $G$ is isomorphic to the Tate module of an abelian variety over $\Qp$ if and only if its associated filtered module satisfies conditions $(1)$---$(4)$ of theorem~\ref{mainthm}.
\end{cor}

\begin{proof}
According to~\cite{Se-Ta} $\S$2 Cor.2(a) the action of inertia is tame when $p>2d +1$. 
\end{proof}

\begin{rem}
Let $V_{\l}$ be an $\l$-adic representation of $G$, $\l\neq p$, with good reduction over the finite Galois extension $K$. It is completely determined by its associated $\l$-adic Weil representation $\Delta_{\l}=\Hom_{\Q_{\l}[I_K]}(V_{\l},\Q_{\l})$ (\cite{Fo3}, contravariant version). Then $V_{\l}$ is isomorphic to the $\l$-adic Tate module of an abelian variety over $\Qp$ if and only if there exists a weakly admissible filtered $\bigl(\varphi,\gal(K/\Qp)\bigr)$-module $(D,\Fil)$ satisfying the conditions of theorem~\ref{mainthm} such that $\W(D)=\Delta_p$ and $\Delta_{\l}$ are compatible. This follows from the compatibility of the system $\bigl(V_{\l}({\mathcal A}_0)\bigr)_{\l}$ where $\l$ runs over all the primes. Obviously this criterion is not handy. However the results in~\cite{No} give some hints in this direction, under the assumption that the ramification degree of $K$ is less than $p-1$ (in particular it is tame and $p\neq 2$). 
\end{rem}

\bigskip

\end{document}